\documentclass{amsart}

\usepackage{amssymb}
\usepackage{latexsym}
\usepackage{amsmath}
\usepackage{euscript}
\usepackage{graphics}

\def\sJ{{\mathfrak J}}

      \def\dC{{\mathbb C}}

   \def\dN{{\mathbb N}}   
      \def\dR{{\mathbb R}}

   \def\dZ{{\mathbb Z}}

\def\cJ{{\mathcal J}}      
   \def\cN{{\mathcal N}}   
\def\cP{{\mathcal P}}

\def\wt#1{{{\widetilde #1} }}
\def\wh#1{{{\widehat #1} }}

\def\bm\chi{\mbox{\boldmath$\chi$}}

\def\e{{\rm e}}
\let\xker=\ker \def\ker{{\xker\,}}

\def\supp{{\rm supp\,}}

\def\e{\varepsilon}
\def\wB{\widetilde{B}}

\def\wb{\widetilde{b}}

\def\kk{\kappa}
\def\f{\varphi}

\newcommand{\sddots}{\begin{picture}(2,2)
\multiput(0,0)(1.5,1){3}{.}\end{picture}}

\def\nk{{\mathbf N}_\kappa}

\def\deg{\operatorname{deg}}

\def\dist{\operatorname{dist}}
\newcommand{\adots}{\makebox[0.9ex][l]{\raisebox{-0.2ex}{.}} \raisebox{0.4ex}{.}
\makebox[0.9ex][r]{\raisebox{1.1ex}{.}}}

\newtheorem{thm}{Theorem}[section]
\newtheorem{cor}[thm]{Corollary}
\newtheorem{lem}[thm]{Lemma}
\newtheorem{prop}[thm]{Proposition}
\theoremstyle{definition}
\newtheorem{defn}[thm]{Definition}
\theoremstyle{remark}
\newtheorem{rem}[thm]{Remark}
\newtheorem{ex}{Example} 
\numberwithin{equation}{section}

\theoremstyle{definition}

\numberwithin{equation}{section}

\begin{document}

\title[Convergence of Pad\'e approximants for definitizable functions ]
{Convergence of diagonal Pad\'e approximants for a class of
definitizable functions }

\author{Maxim Derevyagin}
\address{
Institute of Applied Mathematics and Mechanics \\
R. Luxemburg str. 74\\
83114 Donetsk\\ Ukraine} \email{derevyagin.m@gmail.com}
\author[Vladimir Derkach]{Vladimir Derkach}
\address{Department of Mathematics \\
Donetsk National University \\
Universitetskaya str. 24 \\
83055 Donetsk\\ Ukraine}  \email{derkach.v@gmail.com}

\dedicatory{\it Dedicated to the memory of  Peter Jonas}
\thanks{The research was supported by the Deutschen Akademischen Austauschdienst
 and the  Technical University of Berlin.}
 \subjclass{Primary 30E05; Secondary 47A57.}
\keywords{ Generalized Nevanlinna function, Definitizable function,
Pad\'e approximant, Pontryagin space, Jacobi matrix, Orthogonal
polynomials.}
\begin{abstract}
Convergence of diagonal Pad\'e approximants is studied for a class
of functions which admit the integral representation $ {\mathfrak
F}(\lambda)=r_1(\lambda)\int_{-1}^1\frac{td\sigma(t)}{t-\lambda}+r_2(\lambda),
$ where $\sigma$ is a finite nonnegative measure on $[-1,1]$, $r_1$,
$r_2$ are real rational functions bounded at $\infty$, and $r_1$ is
nonnegative for real $\lambda$. Sufficient conditions for the
convergence of a subsequence of diagonal Pad\'e approximants of $
{\mathfrak F}$ on $\dR\setminus[-1,1]$ are found. Moreover, in the
case when $r_1\equiv 1$, $r_2\equiv 0$ and $\sigma$ has a gap
$(\alpha,\beta)$ containing 0, it turns out that this subsequence
converges in the gap.  The proofs are based on the operator
representation of diagonal Pad\'e approximants of $ {\mathfrak F}$
in terms of the so-called generalized Jacobi matrix associated with
the asymptotic expansion of $ {\mathfrak F}$ at infinity.
\end{abstract}
\maketitle

\section{Introduction}

Let
$F(\lambda)=\displaystyle{-\sum\limits_{j=0}^{\infty}\frac{s_{j}}{{\lambda}^{j+1}}}$
be a formal power series with $s_j\in\dR$, and let $L$, $M$ be
positive integers.  An $[L/M]$ Pad\'e approximant for $F$ is defined
as a ratio
\[
F^{[L/M]}(\lambda)=\frac{A^{[L/M]}\left(\frac{1}{\lambda}\right)}
{B^{[L/M]}\left(\frac{1}{\lambda}\right)}
\]
of polynomials $A^{[L/M]}$, $B^{[L/M]}$ of formal degree $L$ and
$M$, respectively, such that $B^{[L/M]}(0)\ne 0$ and
\begin{equation}\label{eq:1.1}
    \sum_{j=0}^{L+M-1}\frac{s_{j}}{{\lambda}^{j+1}}+F^{[L/M]}(\lambda)= O\left(\frac{1}{{\lambda}^{L+M+1}}\right),\quad \lambda\to\infty.
\end{equation}

The classical Markov theorem~\cite{Mar48} states that for every
nonnegative measure $\sigma$ on the interval $[-1,1]$ and the
function
\begin{equation}\label{markovf}
    F(\lambda)
    =\int_{-1}^1\frac{d\sigma(t)}{t-\lambda}
    \end{equation}
with the Laurent expansion
$F(\lambda)=\displaystyle{-\sum\limits_{j=0}^{\infty}\frac{s_{j}}{{\lambda}^{j+1}}}$
at $\infty$ the diagonal Pad\'e approximants $F^{[n/n]}$ exists for
every $n\in\dN$ and converge to $F$ locally uniformly on
$\dC\setminus[-1,1]$. However, it should be noted  that in the case
when $\sigma$ has a gap $(\alpha,\beta)$ in its support, the
diagonal Pad\'e approximants
$F^{[n/n]}$ do not usually converge
inside the gap (see~\cite{Mar48}).

In~\cite{BakerGW} it was conjectured that for every function $F$
holomorphic in a neighborhood of $\infty$ there is a subsequence of
diagonal $[n/n]$ Pad\'e approximants which converges to $F$ locally
uniformly in the neighborhood of $\infty$ (Pad\'e hypothesis). In
general, as was shown by D.~Lubinsky~\cite{Lub03} (see also~\cite{Bus02}), this conjecture
fails to hold,  but for some classes of functions the Pad\'e
hypothesis is still true. For example, if $F$ has the form
\[
F(\lambda)=\int_{-1}^1\frac{d\sigma(t)}{t-\lambda}+r(\lambda),
\]
where $r$ is a rational function with poles outside of $[-1,1]$, the
convergence of Pad\'e approximants  was proved by
A.~Gonchar~\cite{Gon} and E.~Rakhmanov~\cite{Rak}.

In~\cite{DD,DD07} we studied the Pad\'e hypothesis in the class
of generalized Nevanlinna functions introduced in~\cite{KL79} (see
the definition at the beginning of Section~2), which contains, in particular,
functions of the form
\begin{equation}\label{pertf}
F(\lambda)=r_1(\lambda)\int_{-1}^1\frac{d\sigma(t)}{t-\lambda}+r_2(\lambda),
\end{equation}
where:
\begin{enumerate}
    \item[(A1)] $\sigma$ is a finite nonnegative measure on $[-1,1]$;
    \item[(A2)] $r_1=q_1/w_1$
 is a rational function, nonnegative for real
$\lambda$ ($\deg q_1\le\deg w_1$);
    \item[(A3)] $r_2=q_2/w_2$  is a real
rational function such that $\deg q_2<\deg w_2$.
\end{enumerate}
Let $F$ have the Laurent expansion $
F(\lambda)=\displaystyle{-\sum\limits_{j=0}^{\infty}\frac{s_{j}}{{\lambda}^{j+1}}}
$ at $\infty$, and let $\cN({\bf s})$ 
be the set of all \textit{normal indices} of the sequence ${\bf
s}=\{s_i\}_{i=0}^\infty$, i.e. natural numbers
$n_1<n_2<\dots<n_j<\dots$, for which
\begin{equation}\label{NormInd}
    \det (s_{i+k})_{i,k=0}^{n_j-1} \ne 0, \quad j=1,2,\dots.
\end{equation}
As is known the sequence $\{n_j\}_{j=1}^\infty$ contains all natural
$n$ big enough. The Pad\'e approximants for $F$ were considered
in~\cite{DD07} in connection with the theory of difference equations
\begin{equation}\label{DifEq}
\wt b_{j-1}u_{j-1}-p_j(\lambda)u_{j}+b_{j}u_{j+1}=0,\,\,\,j\in\dN,
\end{equation}
naturally related to the function $F$, where $p_j$ are monic
polynomials of degree $k_j=n_{j+1}-n_j$, $b_j>0$, $\wt b_{j}$ are
real numbers, such that $|\wt b_{j}|=b_j$,
$j\in\dZ_+:=\dN\cup\{0\}$. It turned out that the diagonal Pad\'e
approximants for $F$ exist for all $n=n_j$ and are calculated by the
formula
\[
F^{[n_j/n_j]}(\lambda)=-\frac{Q_{j}(\lambda)}{P_{j}(\lambda)},
\]
where $P_j$, $Q_j$ are polynomials of the first and the second type
associated with the difference equation~\eqref{DifEq}
(see~\cite{DD}). In~\cite{DD07} it was shown that the sequence of
diagonal Pad\'e approximants for $F$ converges to $F$ locally
uniformly in $\dC\setminus ([-1,1]\cup\cP(F))$, where $\cP(F)$ is
the set of poles of $F$. Subdiagonal Pad\'e approximants
$F^{[n_j/n_j-1]}$ for $F$ exist if and only if
\[
n_j\in\cN_F:=\{n_j\in \cN({\bf s}):\,P_{j-1}(0)\ne 0\}.
\]
The convergence  of the sequence $\{F^{[n_j/n_j-1]}\}_{n_j\in
\cN_F}$ of subdiagonal Pad\'e approximants for $F$ on $\dC\setminus
\dR$ was also proved.

In Theorem~\ref{th_subdA} we improve the result
of~\cite[Theorem~4.16]{DD07} by pointing out sufficient condition
for convergence of the sequence of subdiagonal Pad\'e approximants
$F^{[n/n-1]}$ to $F$ in a
neighborhood of $\infty$. 
In the
previous notations this condition takes the form
\begin{enumerate}
    \item[(B)] The sequence $\left\{b_{j-1}\frac{P_j(0)}{P_{j-1}(0)}\right\}_{n_j\in \cN_F}$ is
    bounded.
\end{enumerate}


 The main part of the present paper is dedicated to the study of
convergence of diagonal Pad\'e approximants of definitizable
functions with one "turning point". This class was introduced  by
P.~Jonas in~\cite{J2000}. We postpone the exact definitions until
Section 4 and we mention only that the typical representative of
this class is the function
\begin{equation}\label{pertf2}
{\mathfrak
F}(\lambda)=r_1(\lambda)\int_{-1}^1\frac{td\sigma(t)}{t-\lambda}+r_2(\lambda),
\end{equation}
where $\sigma$, $r_1$, $r_2$  satisfy the assumptions (A1), (A2) and
\begin{enumerate}
    \item[(A$3'$)] $r_2=q_2/w_2$  is a real
rational function such that $\deg q_2\le\deg w_2$.
\end{enumerate}
We prove that for every ${\mathfrak F}$ of the form~\eqref{pertf2}
satisfying condition~(B) the Pad\'e hypothesis is still true. The idea is that the diagonal
Pad\'e approximants for ${\mathfrak F}$ are proportional to the
subdiagonal Pad\'e approximants for $ F$. This fact was observed by
A.~Magnus in a formal setting~\cite{Mag1} and its operator
interpretation 
for the Nevanlinna class was given by B.~Simon~\cite{Si}. This
observation and~\cite[Theorem~4.16]{DD07} enable us to prove in
Theorem~\ref{PAfordef} that the sequence of diagonal Pad\'e
approximants for ${\mathfrak F}$ converges to ${\mathfrak F}$
locally uniformly in $\dC\setminus (\dR\cup\cP({\mathfrak F}))$.
Moreover, using the result of Theorem~\ref{th_subdA} we show that
condition~(B) is sufficient for the convergence of diagonal Pad\'e
approximants of ${\mathfrak F}$ in a neighborhood of $\infty$
(see Theorem \ref{theoremA}).

In Theorem~\ref{theoremC} we specify this result to the case when
 the function $\mathfrak F$ in~\eqref{pertf} takes
 the form
\begin{equation}\label{eq:1.5}
    {\mathfrak F}(\lambda)=\int_E\frac{td\sigma(t)}{t-\lambda},
    \quad E=[-1,\alpha]\cup[\beta,1],
\end{equation}
where the measure $\sigma$ has a gap $ (\alpha,\beta)$
with $\alpha<0<\beta$.
For this function one can observe a new effect, that the sequence of
diagonal Pad\'e approximants \linebreak
$\{{\mathfrak F}^{[n_j-1/n_j-1]}\}_{n_j\in\cN_F}$
converges to ${\mathfrak F}$ in
the gap $(\alpha,\beta)$. The proof of this result is based on the
theory of generalized Jacobi matrices associated with generalized
Nevanlinna functions and on the operator representation of the
subdiagonal Pad\'e approximants for the generalized Nevanlinna
function $F(\lambda):={\mathfrak F}(\lambda)/\lambda$. Moreover, in
Theorem~\ref{theoremC} we prove that for such a function
condition~(B) is also necessary and sufficient for the convergence
of the sequence $\{{\mathfrak F}^{[n_j-1/n_j-1]}\}_{n_j\in\cN_F}$ to
${\mathfrak F}$ in a neighborhood of $\infty$.

This theorem makes a bridge to the theory of classical orthogonal
polynomials. In Proposition~\ref{P_ac} we show that the
condition~(B) 
is in force, if $0$ is not an accumulation point of zeros of
polynomials $P_n$ orthogonal with respect to $\sigma$. In the case
when the measure $\sigma$ in \eqref{eq:1.5} satisfies the Szeg\"o
condition on each of the intervals $[-1,\alpha]$ and $[\beta,1]$ we
inspect the question: under what conditions $0$ is not an
accumulation point of zeros of polynomials $P_n$? In
Proposition~\ref{prop:5.9} we show that the results of
E.~Rakhmanov~\cite{Rak} can be applied to give a partial answer to
this question and, hence, to find some sufficient conditions on
$\alpha$, $\beta$ and $\sigma$ for the existence of a subsequence of
diagonal Pad\'e approximants ${\mathfrak F}^{[n/n]}$ which converges
to $F$ locally uniformly in a neighborhood of $\infty$.

The paper is organized as follows. In Section \ref{sect2} the basic
facts concerning generalized Nevanlinna functions and their operator
representations in terms of generalized Jacobi matrices are given.
In Section \ref{sect3} we state and improve some results
from~\cite{DD07} on the locally uniform convergence of subdiagonal
Pad\'e approximants for generalized Nevanlinna functions. In
Section~\ref{sect4} we introduce the class ${\bf D}_{\kappa,
-\infty}$ of definitizable functions with one "turning point", and
find the formula connecting diagonal Pad\'e approximants for
${\mathfrak F}\in{\bf D}_{\kappa, -\infty}$ with subdiagonal Pad\'e
approximants for generalized Nevanlinna function
$F(\lambda)={\mathfrak F}(\lambda)/\lambda$. In Section \ref{sect5}
we apply our results to subclasses of definitizable functions of the
form~\eqref{pertf2} and \eqref{eq:1.5}. 

This paper is dedicated to the memory of Peter Jonas. Discussions
with him during several of our visits to Berlin have had a
significant influence on the development of this paper.

\section{Preliminaries}
\label{sect2}
\subsection{Moment problem in the class of generalized Nevanlinna functions.}
Let $\kappa$ be a nonnegative integer. 
Recall that a
function $F$, meromorphic in $\dC_+\cup\dC_-$, is said to belong to
the class ${\mathbf N}_\kappa$ if the domain of holomorphy $\rho(F)$
of the function $F$ is symmetric with respect to $\dR$,
$F(\bar\lambda)=\overline{F(\lambda)}$ for $\lambda\in\rho(F)$, and
the kernel
\[
\begin{cases}
       {\sf N}_{F}(\lambda,\omega)=
\frac{F(\lambda)-\overline{F(\omega)}}{\lambda-\overline{\omega}},& \text{$\lambda,\omega\in\rho(F)$};\\
  {\sf N}_{F}(\lambda,\overline{\lambda})=F'(\lambda),& \text{$\lambda\in\rho(F)$}\\
  \end{cases}
\]
has  $\kappa$ negative squares on  $\rho(F)$. The last statement
means that for every $n\in\dN$ and
$\lambda_{1},\lambda_{2},\dots,\lambda_{n}\in\rho(F)$, $n\times n$
matrix $({\sf N}_{F}(\lambda_{i},\lambda_{j}))_{i,j=1}^n$ has at
most $\kappa $ negative eigenvalues (with account of
multiplicities) and for some choice of $n$,
$\lambda_{1},\lambda_{2},\dots,\lambda_{n}$ it has exactly $\kappa$
negative eigenvalues (see~\cite{KL79}).

We will say (cf.~\cite{DHS01}) that a generalized Nevanlinna
function $F$ belongs to
the class ${\bf N}_{\kk,-2n}$ if $F\in\nk$ and for some numbers $
s_0,\dots,s_{2n}\in\dR$ the
following asymptotic expansion holds true
\begin{equation}\label{asymp}
F(\lambda)= 
-\frac{s_{0}}{\lambda}-\frac{s_{1}}{\lambda^2}
-\dots-\frac{s_{2n}}{\lambda^{2n+1}}+o\left(\frac{1}{\lambda^{2n+1}}\right),
\quad\lambda\widehat{\rightarrow }\infty,
\end{equation}
where $\lambda\widehat{\rightarrow }\infty$ means
that $\lambda$ tends to $\infty$ nontangentially, that is inside the
sector $\varepsilon<\arg \lambda<\pi-\varepsilon$ for some
$\varepsilon>0$. Let us set
\[
{\bf N}_{\kk,-\infty}:=\bigcap_{n\ge 0}{\bf N}_{\kk,-2n}.
\]
In particular, every function of the form~\eqref{pertf} where $r_1$,
$r_2$, $\sigma$ are subject to the assumptions (A1)--(A3), belongs
to the class ${\bf N}_{\kk,-\infty}$  for some $\kappa\in \dZ_+$
(see~\cite{KL79}). Moreover, every generalized Nevanlinna function
$F\in{\bf N}_{\kk,-\infty}$ holomorphic at $\infty$ admits the
representation~\eqref{pertf}  for some $r_1$, $r_2$, $\sigma$
satisfying (A1)--(A3).

It will be sometimes convenient to use the following notation
\begin{equation}\label{eq:asymp}
F(\lambda)\sim-\sum_{j=0}^\infty
\frac{s_j}{{\lambda}^{j+1}},\quad\lambda\widehat{\rightarrow }\infty
\end{equation}
to denote the validity of~\eqref{asymp} for all $n\in\dN$.

 \noindent {\bf Problem} ${\bf
M}_{\kk}({\bf s})$. {\it Given are a nonnegative integer $\kk$ and a
sequence ${\bf s}= \{s_j\}_{j=0}^{\infty}$ of real numbers, such
that the matrices $S_{n}:=(s_{i+j})_{i,j=0}^{n}$ are nondegenerate
for all $n$ big enough. Find a function $F\in{\bf N}_{\kk}$, which
has the asymptotic expansion~\eqref{eq:asymp}.}

We say that the problem ${\bf M}_{\kk}({\bf s})$ is {\it
determinate} if  ${\bf M}_{\kk}({\bf s})$ has a unique solution. The
problem ${\bf M}_{\kk}({\bf s})$ was considered in~\cite{KL81},
where it was shown that the problem ${\bf M}_{\kk}({\bf s})$ is
solvable if and only if the number of negative eigenvalues of the
matrix $S_{n}$ does not exceed $\kappa$ for all $n\in\dN$. The Schur
algorithm for solving the problem ${\bf M}_{\kk}({\bf s})$
considered in~\cite{De} proceeds as follows.

Let $\cN({\bf s})$ 
be the set of all \textit{normal indices} of the sequence ${\bf s}$,
i.e. natural numbers $n_j$, for which
\begin{equation}\label{NormInd1}
    \det S_{n_j-1} \ne 0,\quad j=1,2,\dots.
\end{equation}
If $F_0:=F$ is a generalized Nevanlinna function with the asymptotic
expansion~\eqref{eq:asymp}, then the function $-1/F_0$ can be
represented as
\begin{equation}\label{SchurTr}
-\frac{1}{F_0(\lambda)}=\varepsilon_0p_0(\lambda)+b_0^2 F_1(\lambda),
\end{equation}
where $\varepsilon_0=\pm 1$, $b_0>0$, $p_0$ is a monic polynomial of
degree $k_0=n_1$ and $F_1$ is a generalized Nevanlinna function.
Continuing this process one gets sequences $\varepsilon_j=\pm 1$,
$b_j>0$, $j\in\dZ_+$ and a sequence of real monic
polynomials $p_j$ of degree $k_j=n_{j+1}-n_j$, such that $F$ admits
the following expansion into a $P$-fraction
\begin{equation}\label{Pfraction}
-\frac{\varepsilon_0}{p_0(\lambda)}
\begin{array}{l} \\ -\end{array}
\frac{\varepsilon_0\varepsilon_1b_0^2}{p_1(\lambda)}
\begin{array}{ccc} \\ - & \cdots & -\end{array}
\frac{\varepsilon_{N-1}\varepsilon_N
b_{N-1}^2}{p_N(\lambda)}\begin{array}{l}\\-\dots\end{array}.
\end{equation}
 A similar algorithm for a continued fraction expansion of a formal power
series was proposed by A.~Magnus in~\cite{Mag1}. The objects
$\varepsilon_j$, $b_j$, $p_j$ are uniquely defined by the sequence
${\bf s}= \{s_j\}_{j=0}^{\infty}$ (see~\cite{De}, \cite{DD07}). The
function $F_1$ in~\eqref{SchurTr} is called the Schur transform of
$F_0\in {\bf N}_\kappa$ (cf.~\cite{ADL07}).


\subsection{Generalized Jacobi matrices}
 Let
$p(\lambda)=p_{k}{\lambda}^{k}+\dots+p_{1}\lambda+p_{0}$, $p_{k}=1$ be a monic scalar
real polynomial of degree $k$. Let us associate with the polynomial $p$ its
symmetrizator $E_p$ and the companion matrix $C_p$, given by
\begin{equation}\label{comp}
E_{p}=\begin{pmatrix}
p_{1}&\dots&p_{k}\\
\vdots&\sddots&\\
p_{k}&&\\
\end{pmatrix},\quad
C_{p}=\begin{pmatrix}
0   &\dots  &0      &-p_{0}\\
1   &       &        &-p_{1}\\
    &\ddots &      &\vdots\\
&&1&-p_{k-1}\\
\end{pmatrix},
\end{equation}
where all the non-specified elements are supposed to be equal to
zero. As is known (see~\cite{GLR}), $\det(\lambda-C_p)=p(\lambda)$
and the spectrum $\sigma(C_p)$ of the companion matrix $C_p$ is
simple. The matrices $E_{p}$ and $C_{p}$ are related by
\begin{equation}\label{cb}
C_{p}E_{p}=E_{p}C_{p}^{\top}.
\end{equation}
\begin{defn}
Let $p_j$ be real monic polynomials of degree ${k_j}$
\[
p_{j}(\lambda)={\lambda}^{k_{j}}+p_{k_{j}-1}^{(j)}{\lambda}^{k_{j}-1}\dots+
p_{1}^{(j)}\lambda+p_{0}^{(j)},
\]
and let $\e_{j}=\pm 1$, $b_{j}>0$, $j\in\dZ_+$. The tridiagonal
block matrix
\begin{equation}\label{mJacobi}
\cJ=\begin{pmatrix}
A_{0}   &\wB_{0}&       &\\
B_{0}   &A_1    &\wB_{1}&\\
        &B_1    &A_{2} &\ddots\\
&       &\ddots &\ddots\\
\end{pmatrix}
\end{equation}
where $A_{j}=C_{p_{j}}$ and $k_{j+1}\times k_{j}$ matrices $B_{j}$
and $k_{j}\times k_{j+1}$ matrices $\wB_{j}$ are given by
\begin{equation}\label{bblock}
B_{j}=\begin{pmatrix}
0&\dots&b_{j}\\
\hdotsfor{3}\\
0&\dots&0\\
\end{pmatrix},\,
\wB_{j}= \begin{pmatrix}
0&\dots&\wt b_{j}\\
\hdotsfor{3}\\
0&\dots&0\\
\end{pmatrix}, \,
\wt b_{j}=\e_j\e_{j+1} b_j, \, j\in\dZ_+.
\end{equation}
will be called a {\it generalized Jacobi matrix}  associated with
the sequences of polynomials  $\{{\e_{j}p_{j}\}}_{j=0}^{\infty}$
and numbers $\{b_j\}_{j=0}^{\infty}$.
\end{defn}
\begin{rem} Define an infinite matrix $G$ by the
equality
\begin{equation}\label{Gram}
G=\mbox{diag}(G_{0},\dots,G_{N},\dots),\quad
G_{j}=\e_{j}E_{p_{j}}^{-1},\,\, j=0,\dots, N
\end{equation}
and let  $\ell^2_{[0,\infty)}(G)$ be the space of $\ell^2$-vectors
with the inner product
\begin{equation}\label{metric}
\left[x,y\right]=(Gx,y)_{\ell^2_{[0,\infty)}},\quad x,y\in\ell^2_{[0,\infty)}.
\end{equation}
The inner product~\eqref{metric} is indefinite, if either $k_j>1$
for some $j\in\dZ_+$, or at least one $\e_{j}$ is equal to $-1$. The
space $\ell^2_{[0,\infty)}(G)$ is equivalent to a Kre\u{\i}n space
(see~\cite{AI}) if both $G$ and $G^{-1}$ are bounded in
$\ell^2_{[0,\infty)}$.  If $k_j=\e_j=1$ for all $j$ big enough, then
$\ell^2_{[0,\infty)}(G)$ is a Pontryagin space. In these cases it
follows from~\eqref{cb} that the generalized Jacobi matrix $\cJ$
determines a symmetric operator $S$ in the space
$\ell^2_{[0,\infty)}(G)$ (see details
in~\cite{DD}). 
More general definition of generalized Jacobi matrix which is not connected with the
Schur algorithm has been considered in~\cite{KL79}.
\end{rem}

Setting $\wt b_{-1}=\varepsilon_0$, define polynomials of the
first kind (cf.~\cite{KL79}) $P_j(\lambda)$, $j\in\dZ_+$, as solutions
$u_j=P_j(\lambda)$ of the following system:
\begin{equation}\label{eq10}
\wt b_{j-1}u_{j-1}-p_j(\lambda)u_{j}+b_{j}u_{j+1}=0,\,\,\,j\in\dZ_+,
\end{equation}
with the initial conditions
\begin{equation}\label{InConP}
u_{-1}=0,\quad u_0=1.
\end{equation}
Similarly, the polynomials of the second kind $Q_j(\lambda)$,
$j\in\dZ_+$, are defined as solutions $u_j=Q_j(\lambda)$ of the
system~\eqref{eq10} subject to the following initial conditions
\begin{equation}\label{InConP2}
u_{-1}=-1,\quad u_0=0.
\end{equation}
It follows from~\eqref{eq10} that $P_{j}$ is a polynomial of degree
$n_j=\sum\limits_{i=0}^{j-1}k_{i}$ with the leading coefficient $(b_0\dots b_{j-1})^{-1}$
and $Q_j$ is a polynomial of degree $n_j-{k_0}$ with the leading
coefficient $\e_0(b_0\dots b_{j-1})^{-1}$. 
The equations~\eqref{eq10}
coincide with the three-term recurrence relations associated with $P$-fractions
(\cite{Mag1}, see also~\cite[Section 5.2]{JTh}). The following statement is immediate
from~\eqref{eq10}.
\begin{prop}\label{Interlace} {\rm (\cite{DD07}).}
    Polynomials $P_{j}$ and $P_{j+1}$ ($Q_{j}$ and $Q_{j+1}$) have no common zeros.
\end{prop}

 The following connection between the polynomials of the
first and second kinds $P_{j}$, $Q_{j}$ and the shortened Jacobi matrices
$\cJ_{[0,j]}$ can be found in~\cite[Proposition 3.3]{DD07} (in the classical case
see~\cite[Section 7.1.2]{Be}).
\begin{prop}\label{detf}
Polynomials $P_{j}$ and $Q_j$ can be found by the formulas
\begin{eqnarray}
P_{j}(\lambda)&=&(b_{0}\dots b_{j-1})^{-1}\det(\lambda-\cJ_{[0,j-1]}),\label{polynom1}\\
Q_{j}(\lambda)&=&\e_0(b_0\dots b_{j-1})^{-1}\det(\lambda-\cJ_{[1,j-1]}).\label{polynom2}
\end{eqnarray}
\end{prop}
Clearly, $\cJ_{[0,j-1]}$ is a symmetric operator in the subspace $\ell^2_{[0,j-1]}(G)$
of the indefinite inner product space $\ell^2_{[0,\infty)}(G)$, which consists of
vectors $u=\{{u_{ik}}\}_{i=0,\dots,\infty}^{k=0,\dots,n_i-1}$ such that $u_{ik}=0$ for
$i\ge j$. The m-function of the shortened matrix $\cJ_{[0,j-1]}$ is defined by
\begin{equation}\label{SWfun}
    m_{[0,j-1]}(\lambda)=[(\cJ_{[0,j-1]}-\lambda)^{-1}e,e].
\end{equation}
Due to formulas~\eqref{polynom1}, \eqref{polynom2} it is
calculated by
\begin{equation}~\label{mQP}
m_{[0,j-1]}(\lambda)=
-\varepsilon_0\frac{\det(\lambda-\cJ_{[1,j-1]})}{\det(\lambda-\cJ_{[0,j-1]})}=-\frac{Q_{j}(\lambda)}{P_{j}(\lambda)}.
\end{equation}
\begin{rem}\label{ResSet}
    Let us emphasize that the polynomials $P_{j}$ and $Q_{j}$ have no common
    zeros (see~\cite[Proposition 2.7]{DD})
    and due to~\eqref{mQP} the set of
    holomorphy of $m_{[0,j-1]}$ coincides with the resolvent set of $\cJ_{[0,j-1]}$.
\end{rem}

\begin{thm}[\cite{DD07}]\label{DiagPA}
Let $F\in{\bf N}_{\kappa,-\infty}$ and  the corresponding indefinite
moment problem $M_{\kk}({\bf s})$ be determinate. Then:
\begin{enumerate}
    \item[(i)] the
generalized Jacobi matrix corresponding to $F$ via~\eqref{Pfraction}
and~\eqref{mJacobi} generates a selfadjoint operator $\cJ$ in
$\ell^2_{[0,\infty)}(G)$ and
\begin{equation}\label{Wfun}
    F(\lambda)=[(\cJ-\lambda)^{-1}e,e].
\end{equation}
    \item[(ii)] the diagonal $[n_j/n_j]$ Pad\'e
approximants of $F(\lambda)$ coincide with $m_{[0,j-1]}(\lambda)$ and
converge to $F(\lambda)$ locally uniformly on $\dC\setminus\dR$.
\end{enumerate}
\end{thm}
The proof of this result is based on the fact that the compressed
resolvents of $\cJ_{[0,j-1]}$ converge to the compressed resolvent
of $\cJ$ (see~\cite[Theorem 4.8]{DD07}). Theorem~\ref{DiagPA}
contains as partial cases some results of A.~Gonchar~\cite{Gon} and\linebreak
E.~Rakhmanov~\cite{Rak}, mentioned  in Introduction, as well as the
results of \linebreak G.L.~Lopes~\cite{Lop} concerning convergence of diagonal
Pade approximants for rational perturbations of Stieltjes functions.

\section{The  convergence of subdiagonal Pad\'e approximants.}
\label{sect3}
Let us consider the following finite generalized Jacobi matrix
\begin{equation}\label{pert_jacobi}
\cJ_{[0,j]}(\tau)=
\begin{pmatrix}
A_{0}&\wB_{0}&&\\
B_{0}&\ddots&\ddots&\\
&\ddots&A_{j-1}&\wB_{j-1}\\
&&B_{j-1}&A_{j}(\tau)\\
\end{pmatrix},\,\,
A_j(\tau)=\begin{pmatrix}
0&\dots&0&-p_{0}^{(j)}+\tau\\
1&&&-p_{1}^{(j)}\\
&\ddots&&\vdots\\
&&1&-p_{k_j-1}^{(j)}\\
\end{pmatrix}.
\end{equation}
A vector $u=(u_{ik})_{i=0,\dots,j}^{k=0,\dots,n_i-1}\in\dC^{n_j}$ is
a left eigenvector of the matrix $\cJ_{[0,j]}(\tau)$, corresponding
to the eigenvalue 0 if and only if  $u_i=u_{i0}$, $i=0,\dots,j$,
satisfy the equations~\eqref{eq10} for $i=0,\dots,j-1$, and
\begin{equation}\label{zero_equation}
\wb_{j-1}u_{j-1}-(p_j(0)-\tau)u_j=0,
\end{equation}
\[
u_{ik}=0,\quad i=0,\dots,j,\,\, k=1,\dots,n_i-1.
\]
Therefore $u_{i}=P_{i}(0)$ $(i=0,\dots,j)$. If $u_j=P_j(0)\ne 0$,
then~\eqref{zero_equation} yields
\begin{equation}\label{taurepr}
\tau=\tau_j:=-\frac{\wb_{j-1}P_{j-1}(0)-p_{j}(0)P_j(0)}{P_j(0)}=
b_j\frac{P_{j+1}(0)}{P_j(0)}.
\end{equation}
\begin{prop}\label{zero_jacobi}
If $P_j(0)\ne 0$ then there exists a number $\tau_j\in\dR$ such that
\begin{equation}\label{0sp}
0\in\sigma_p(\cJ_{[0,j]}(\tau_j)).
\end{equation}
Moreover, $\tau_j$ can be found by the formula~\eqref{taurepr}.
\end{prop}

\begin{rem}\label{zero_jacobi2}
We will use the notation $\cJ_{[0,j]}^{(K)}$ for the matrix
$\cJ_{[0,j]}(\tau_j)$ with the property~\eqref{0sp}. In the case
when the corresponding indefinite moment problem $M_{\kk}({\bf s})$
is determinate for all $\lambda\in\dC\setminus(\dR\cup\cP(F))$ we
have
\begin{equation}\label{SWfunK}
  m_{[0,j-1]}^{(K)}(\lambda):=[(\cJ_{[0,j-1]}^{(K)}-\lambda)^{-1}e,e]\to F(\lambda)=[(\cJ-\lambda)^{-1}e,e]
\end{equation}
as $j\to\infty$ (see~\cite[Proposition 4.4]{DD07}).
\end{rem}
\begin{thm}[\cite{DD07}]\label{t_subdiagonal}
Let a function $F\in {\bf N}_\kappa$ have the expansion~\eqref{eq:asymp}
for every $n\in\dN$, let the corresponding moment problem ${\bf
M}_{\kappa}({\bf s})$ be determinate, and let
$\{n_j\}_{j=1}^\infty=\cN({\bf s})$ be the set of normal indices of
the sequence ${\bf s}=\{s_{i}\}_{i=0}^\infty$. Then:
\begin{enumerate}
\item[(i)] The $[n_j/n_j-1]$ Pad\'e approximant $F^{[n_j/n_j-1]}$
exists if and only if
\begin{equation}\label{cN}
n_j\in \cN_F:=\{n_{j}\in\cN({\bf s}):\,P_{j-1}(0)\ne 0\};
\end{equation}
\item[(ii)] The sequence
\[
F^{[n_j/n_j-1]}=-\frac{Q_j(\lambda)P_{j-1}(0)-Q_{j-1}(\lambda)P_{j}(0)}
{P_j(\lambda)P_{j-1}(0)-P_{j-1}(\lambda)P_{j}(0)},\quad n_j\in \cN_F,
\]
converges to $F$ locally uniformly in
$\dC\setminus(\dR\cup\cP(F))$.
\end{enumerate}
\end{thm}
\begin{proof}
We will sketch the proof of this theorem presented in~\cite{DD07}.

Let $P_{j-1}(0)\ne 0$ and let $P_j^{(K)}$, $Q_j^{(K)}$ be the
polynomials of the first and second kinds, respectively, associated
with the matrix $\cJ_{[0,j-1]}^{(K)}$, that is
\[
\begin{split}P_{j}^{(K)}(\lambda)&=(b_{0}\dots
b_{j-1})^{-1}\det(\lambda-\cJ_{[0,j-1]}^{(K)}),\\
Q_{j}^{(K)}(\lambda)&=(b_{0}\dots b_{j-1})^{-1}\det(\lambda-\cJ_{[1,j-1]}^{(K)}). \end{split}\]
Next, applying the determinant decomposition theorem one can obtain
\[
\begin{split}
P_{j}^{(K)}(\lambda)&=P_j(\lambda)-\frac{\tau_{j-1}}{b_{j-1}}P_{j-1}(\lambda)
=\frac{P_j(\lambda)P_{j-1}(0)-P_{j-1}(\lambda)P_{j}(0)}{P_{j-1}(0)},\\
Q_{j}^{(K)}(\lambda)&=Q_j(\lambda)-\frac{\tau_{j-1}}{b_{j-1}}Q_{j-1}(\lambda)
=\frac{Q_j(\lambda)P_{j-1}(0)-Q_{j-1}(\lambda)P_{j}(0)}{P_{j-1}(0)}
\end{split}\]
It obviously follows from~\eqref{SWfunK} that
\begin{equation}\label{defmzero}
m_{[0,j-1]}^{(K)}(\lambda):=\left[(\cJ_{[0,j-1]}^{(K)}-\lambda)^{-1}e,e\right]\sim
-\sum_{i=0}^{\infty}\frac{{s_{i}}^{(K)}}{\lambda^{i+1}},\quad \lambda\wh\to\infty,
\end{equation}
where ${s_{i}}^{(K)}=\left[(\cJ_{[0,j-1]}^{(K)})^ie,e\right]$. Using the form of the
matrix $\cJ_{[0,j-1]}^{(K)}$ one gets
\[
s_i^{(K)}=s_i\mbox{ if }i\le 2n_j-2,
\]
\[
s_{2n_j-1}^{(K)}=s_{2n_j-1}+(b_0\dots b_{j-1})^2\varepsilon_{j-1}\tau_{j-1}.
\]
So, the function $m_{[0,j-1]}^{(K)}(\lambda)$ has the following asymptotic expansion
\[
m_{[0,j-1]}^{(K)}(\lambda)= -\sum_{i=0}^{2n_j-2}\frac{s_{i}}{\lambda^{i+1}}+
                O\left(\frac{1}{\lambda^{2n_j}}\right),\quad
                \lambda\wh\to\infty,
\]
where the sequence $\{s_j\}_{j=0}^{\infty}$ corresponds to the generalized Jacobi matrix
$\cJ$. On the other hand, due to~\eqref{mQP}
\[
m_{[0,j-1]}^{(K)}(\lambda)= -
\frac{Q_{j}^{(K)}(\lambda)}{P_{j}^{(K)}(\lambda)}.
\]
Further, setting
\[
A^{[n_i/n_i-1]}\left(\frac{1}{\lambda}\right)=\left(\frac{1}{\lambda}\right)^{n_{i}}Q_{i}^{(K)}(\lambda),\quad
B^{[n_i/n_i-1]}\left(\frac{1}{\lambda}\right)=\left(\frac{1}{\lambda}\right)^{n_{i}}P_{i}^{(K)}(\lambda),
\]
for $i=0,1\dots,j$ and taking into account  the equality $P_{j}^{(K)}(0)=0$, one obtains
\begin{equation}\label{PA_sub}
m_{[0,j-1]}^{(K)}(\lambda)= \frac{A^{[n_{j}/n_{j}-1]}\left({1}/{\lambda}\right)}
{B^{[n_{j}/n_{j}-1]}\left({1}/{\lambda}\right)},
\end{equation}
where
\[
\deg A^{[n_{j}/n_{j}-1]}=n_{j},\quad \deg B^{[n_{j}/n_{j}-1]}=n_{j}-1, \quad
 B^{[n_{j}/n_{j}-1]}(0)=\frac{1}{b_0\dots b_{j-1}}\ne 0.
\]
Therefore,  $m_{[0,j-1]}^{(K)}(\lambda)$ is the $[n_{j}/n_{j}-1]$ Pad\'e
approximant for the corresponding Hamburger series. This proves the
first part of the theorem. The second statement rests on the fact
mentioned in Remark~\ref{zero_jacobi2}.
\end{proof}
\begin{rem}
Condition $P_{j-1}(0)\ne 0$ is equivalent to
\begin{equation}\label{cN2}
    \begin{vmatrix}
s_1 & \dots &  s_{n_{j-1}}\\
\hdotsfor{3}\\
s_{n_{j-1}}& \dots& s_{2n_{j-1}-1}
\end{vmatrix}\ne 0.
\end{equation}
It follows from Proposition~\ref{Interlace} that the set $\cN_F$ is infinite.
\end{rem}
\begin{ex}
Consider the following classical $2$-periodic Jacobi matrix
\[
\cJ=\begin{pmatrix}
a_0       & b_0   &       &    &   \\
b_0       & a_1   & b_1   &    &  \\
          & b_1   & a_2   &  \ddots \\
 &       & \ddots & \ddots        \\
\end{pmatrix},
\quad a_n=\frac{(-1)^n+1}{2},\quad b_n=1, \quad n\in\dZ_+.
\]
The m-function corresponding to this Jacobi matrix can be found by
using standard methods (see~\cite{NikSor})
\[
\f(\lambda)=((\cJ-\lambda)^{-1}e,e)_{\ell^2_{[0,\infty)}}=\frac{\lambda-\lambda^2+\sqrt{(\lambda^2-\lambda-2)^2-4}}{2(\lambda-1)}
\]
(here the branch of the square root is determined by the condition
$\f(\lambda)=-1/\lambda+o\left(1/\lambda\right)$
as $\lambda\to\infty$). Therefore, the function
$\f(\lambda)$ admits the integral representation
\[
    \f(\lambda)=\int_{E}\frac{d\sigma(t)}{t-\lambda}\,,
\]
where the support $E:=\mbox{supp }\sigma$ of the measure $\sigma$ is
contained in $[-2,3]$. Since the $[n/n-1]$ Pad\'e approximant
 is equal to $f^{[n/n-1]}(\lambda)=m_{[0,n-1]}^{(K)}(\lambda)$,
its poles coincide with eigenvalues of the matrix
$\cJ_{[0,n-1]}^{(K)}$. Let us show that the eigenvalue of the matrix
$\cJ_{[0,2k]}^{(K)}$ with the largest absolute value  tends to
infinity as $k\to+\infty$. First we compute $\tau_n$. Since
$\tau_n=b_nP_{n+1}(0)/P_n(0)=-b_{n-1}\wb_{n-1}/\tau_{n-1}+p_n(0)$,
we have
\[
\tau_n=-1/\tau_{n-1}-((-1)^n+1)/2.
\]
Clearly, $\tau_0=-1$. By induction, we have the following formulas
\[
\tau_{2k}=-(k+1),\quad\quad \tau_{2k+1}=1/(k+1).
\]
Taking into account that $\cJ_{[0,2k]}^{(K)}$ is a self-adjoint
matrix, one obtains
\[
|\lambda_{max}(\cJ_{[0,2k]}^{(K)})|=\Vert \cJ_{[0,2k]}^{(K)}\Vert
\ge|(\cJ_{[0,2k]}^{(K)}e_{2k},e_{2k})_{\ell^2}|=k,
\]
where $\lambda_{max}(\cJ_{[0,2k]}^{(K)})$ is the eigenvalue of the
matrix $\cJ_{[0,2k]}^{(K)}$ with the largest absolute value.
Therefore, $|\lambda_{max}(\cJ_{[0,2k]}^{(K)})|\to+\infty$ as
$k\to+\infty$. So, infinity is an accumulation point of the set of
poles of the Pad\'e approximants $f^{[n/n-1]}$ of the function
$\f$ holomorphic at infinity.
(This example was given in~\cite{DD07} with several misprints.)
\end{ex}
Under certain conditions, it is possible to say more about the convergence of the
sequence $F^{[n/n-1]}$ on the real line.

\begin{thm}\label{th_subdA}
Let $F$ have the form
\[
F(\lambda)=r_1(\lambda)\int_{-1}^1\frac{d\sigma(t)}{t-\lambda}+r_2(\lambda)=-\sum_{j=0}^\infty
\frac{s_j}{{\lambda}^{j+1}},\quad |\lambda|>R,
\]
where $\sigma$, $r_1$, $r_2$ satisfy the assumptions (A1)-(A3). If
the sequence $\{\tau_{j-1}\}_{n_j\in\cN_F}$ is bounded, i.e.
\begin{equation}\label{Pcond}
\sup\limits_{n_j\in \cN_F}\left|b_{j-1}\frac{P_{j}(0)}{P_{j-1}(0)}\right|<\infty,
\end{equation}
than there exists a constant $\varepsilon>0$ such that the sequence
$\{F^{[n_j/n_j-1]}\}_{n_j\in \cN_F}$ converges to $F$ locally
uniformly in
$\dC\setminus([-1-\varepsilon,1+\varepsilon]\cup\cP(F))$.
\end{thm}
\begin{proof}
It is obvious that $F$ corresponds to the determinate moment problem ${\bf
M}_{\kappa}({\bf s})$. Moreover, the corresponding generalized Jacobi matrix ${\cJ}$ is
a bounded linear operator. According to Theorem~\ref{t_subdiagonal}, the sequence
$\{F^{[n_j/n_j-1]}\}_{n_j\in \cN_F}$ converges to $F$ locally uniformly in
$\dC\setminus(\dR\cup\cP(F))$. Due to~\eqref{PA_sub} we have
\begin{equation}\label{mK}
F^{[n_j/n_j-1]}(\lambda)=m_{[0,j-1]}^{(K)}(\lambda)=\left[(\cJ_{[0,j-1]}^{(K)}-\lambda)^{-1}e,e\right].
\end{equation}
Since the sequence $\{\tau_{j-1}\}_{n_j\in \cN_F}$ is bounded, one
obtains
\begin{equation}\label{eq:est}
\Vert \cJ_{[0,j-1]}^{(K)}\Vert\le\Vert
\cJ_{[0,j-1]}\Vert+|\tau_{j-1}|\le 1+\varepsilon,\quad n_j\in \cN_F
\end{equation}
for some $\varepsilon>0$. 
It follows from the inequality
\[
\Vert(\cJ_{[0,j-1]}^{(K)}-\lambda)^{-1}\Vert \le\frac{1}
{|\lambda|-\|\cJ_{[0,j-1]}^{(K)}\|}\quad
(|\lambda|>\|\cJ_{[0,j-1]}^{(K)}\|)
\]
and~\eqref{eq:est} that
\[
\begin{split}
\left|m_{[0,j-1]}^{(K)}(\lambda)\right|&=
\left|\left((\cJ_{[0,j-1]}^{(K)}-\lambda)^{-1}e,Ge\right)_{\ell^2}\right|\\
&\le\Vert(\cJ_{[0,j-1]}^{(K)}-\lambda)^{-1}e\Vert_{\ell^2}\Vert
Ge\Vert_{\ell^2} \le \frac{\Vert
Ge\Vert_{\ell^2}}{|\lambda|-1-\varepsilon}.
\end{split}
\]
for $|\lambda|>1+\varepsilon$. To complete the proof, it is
sufficient to apply the Vitali theorem.
\end{proof}

\begin{cor}
If the sequence $\left\{\tau_{j-1}\right\}_{n_j\in \cN_F}$ tends to
$0$ then the sequence \linebreak $\left\{F^{[n_j/n_j-1]}\right\}_{n_j\in
\cN_F}$ converges to $F$ locally uniformly in
$\dC\setminus([-1,1]\cup\cP(F))$.
\end{cor}
\begin{proof}
The statement is implied by~\eqref{SWfun}, \eqref{mK}, the relation
\[
\Vert \cJ_{[0,j-1]}-\cJ_{[0,j-1]}^{(K)}\Vert\to 0 ,\quad j\to\infty,
\]
and Theorem~\ref{DiagPA}.
\end{proof}
\begin{rem}
The condition~\eqref{Pcond} can be reformulated in terms of the monic orthogonal
polynomials
\[
\widehat{P}_j(\lambda)=\frac{1}{\det S_{n_j-1}}\left|%
\begin{array}{cccc}
  s_0       & s_1       & \dots & s_{n_j} \\
 \hdotsfor{4}\\
  s_{n_j-1} & s_{n_j}   & \dots & s_{2n_j-1} \\
  1         & \lambda   & \dots & \lambda^{n_j} \\
\end{array}%
\right|,
\]
which are connected with ${P}_j(\lambda)$ by the formulas $
\widehat{P}_j(\lambda)=(b_0\cdots b_{j-1}){P}_j(\lambda),\, j\in\dN.
$ Therefore, the condition~\eqref{Pcond} takes the form
\begin{equation}\label{WhPcond}
\sup\limits_{n_j\in \cN_F}\left|\frac{\wh P_{j}(0)}{\wh P_{j-1}(0)}\right|<\infty.
\end{equation}
\end{rem}

\begin{rem}\label{Subs}
It is clear from the proof of Theorem~\ref{th_subdA} that the
existence of a converging subsequence of the $[n/n-1]$ Pad\'e
approximants follows from the existence of a bounded subsequence
of $\{\tau_{j-1}\}_{n_j\in \cN_F}$.
\end{rem}

\section{A class of definitizable functions and Pade approximants.}
\label{sect4}
\subsection{Classes ${\bf D}_{\kk,-\infty}$ and ${\bf
D}^\circ_{\kk,-\infty}$.}
\begin{defn}\label{deffunct}
Let us say that a function ${\mathfrak F}$ meromorphic in $\dC_+$ belongs to the class
${\bf D}_{\kk,-\infty}$ if
\[
F(\lambda):=\frac{{\mathfrak F}(\lambda)}{\lambda}\in{\bf N}_{\kk,-\infty}\quad\text{and}\quad
{\mathfrak F}(\lambda)=O(1),\quad\lambda\widehat{\rightarrow }\infty.
\]
\end{defn}
 Clearly, every
function ${\mathfrak F}\in {\bf D}_{\kk,-\infty}$ is definitizable in the sense
of~\cite{J2000}. Indeed, consider the factorization
\[
\frac{{\mathfrak
F}(\lambda)}{\lambda}=r^{-1}(\lambda)(r^{\sharp})^{-1}(\lambda){F}_0(\lambda),
\]
where $r$ is a real rational function,
$r^{\sharp}(\lambda)=\overline{r(\overline{\lambda})}$, and ${F}_0\in {\bf
N}_0$. Then
\[
\frac{r(\lambda)r^{\sharp}(\lambda)}{\lambda}{\mathfrak F}(\lambda)={F}_0(\lambda)\in {\bf
N}_0
\]
and, hence, ${\mathfrak F}$ is definitizable,
${\displaystyle\frac{r(\lambda)r^{\sharp}(\lambda)}{\lambda}}$ is definitizing multiplier.

It follows from~\eqref{asymp} that every function ${\mathfrak F}\in {\bf
D}_{\kk,-\infty}$ admits the asymptotic expansion
\begin{equation}\label{asympF}
{\mathfrak F}(\lambda)\sim-{\mathfrak s}_{-1} -\frac{{\mathfrak s}_{0}}{\lambda}-\frac{{\mathfrak
s}_{1}}{\lambda^{2}}-\dots-\frac{{\mathfrak s}_{2n}}{\lambda^{2n+1}}
-\dots,\quad\lambda\widehat{\rightarrow }\infty
\end{equation}
which is connected with the asymptotic expansion~\eqref{asymp} of
$F(\lambda)=\frac{{\mathfrak F}(\lambda)}{\lambda}$ via the formulas
\begin{equation}\label{sfrak}
   {\mathfrak s}_{j-1}=s_{j},\,\,j\in\dZ_+.
\end{equation}
In what follows we use the Gothic script for all the notations
associated with the ${\bf D}_{\kk,-\infty}$ function and the Roman
script for the ${\bf N}_{\kk,-\infty}$ function to avoid confusion.
We also say that a function ${\mathfrak F}$ meromorphic in $\dC_+$
belongs to the class ${\bf D}^\circ_{\kk,-\infty}$ if
\[
F(\lambda):=\frac{{\mathfrak F}(\lambda)}{\lambda}\in{\bf N}_{\kk,-\infty},\quad
{\mathfrak F}(\lambda)=o(1),\quad\lambda\widehat{\rightarrow }\infty,
\]
and the asymptotic expansion of the function ${\mathfrak F}$
\begin{equation}\label{asympF0}
{\mathfrak F}(\lambda)\sim-\frac{{\mathfrak s}_{0}}{\lambda}-\frac{{\mathfrak
s}_{1}}{\lambda^{2}}-\dots-\frac{{\mathfrak s}_{2n}}{\lambda^{2n+1}}
-\dots,\quad\lambda\widehat{\rightarrow }\infty
\end{equation}
is {\it normalized} in a sense that the first nontrivial coefficient
in~\eqref{asympF0} has modulus 1,
\[
|{\mathfrak s}_{{\mathfrak n}_1-1}|=1.
\]
Let the set $\cN({\mathfrak s})$ of normal indices of the sequence
${\mathfrak s}=\{{\mathfrak s}_{i}\}_{i=0}^\infty$ corresponding to
a function ${\mathfrak F}\in{\bf D}_{\kk,-\infty}^\circ$ be defined
by~\eqref{NormInd1}, that is
\begin{equation}\label{NormIndFrak}
    \cN({\mathfrak s})=\{{\mathfrak n}_j:\det ({\mathfrak s}_{i+k})_{i,k=0}^{{\mathfrak n}_j-1} \ne 0, \quad j=1,2,\dots\}.
\end{equation}
\subsection{Normal indices of the ${\bf D}_{\kk,-\infty}^\circ$ functions.}
Remind that the point $\infty$ is called a generalized pole of
nonpositive type of $F\in{\bf N}_\kappa$ with multiplicity
$\kappa_\infty(F)$, if
\begin{equation}
\label{infgpol}
  0\leq\lim_{\lambda\widehat{\rightarrow }\infty }
   \frac{F(\lambda)}{\lambda^{2\kappa_{\infty}+1}} < \infty,\quad
 -\infty\leq\lim_{\lambda\widehat{\rightarrow }\infty }
   \frac{F(\lambda)}{\lambda^{2\kappa_{\infty}-1}} < 0.
\end{equation}
Similarly, the point $\infty$ is called a generalized zero of
nonpositive type of $F$ with  multiplicity $\pi_\infty(F)$, if
\begin{equation}
\label{infgzer}
  \infty \leq\lim_{\lambda\widehat{\rightarrow }\infty }
   {\lambda^{2\pi_{\infty}+1}}F(\lambda) < 0,\quad
 0\leq \lim_{\lambda\widehat{\rightarrow }\infty }
   {\lambda^{2\pi_{\infty}-1}}F(\lambda) < \infty.
\end{equation}
It was shown in~\cite{KL81} that the multiplicity of $\infty$ as a
generalized pole (zero) of nonpositive type of $F\in {\bf N}_\kappa$
does not exceed $\kappa$.
\begin{lem}\label{lem:4.2}
    Let ${\mathfrak F}\in {\bf D}_{\kk,-\infty}^\circ$,
    let the sequence ${\mathfrak s}=\{{\mathfrak s}_{j}\}_{j=0}^\infty$ be defined
    by the asymptotic expansion~\eqref{asympF0},
    and let $ \cN({\mathfrak s})=\{{\mathfrak n}_j\}_{j=1}^\infty\,$ be the set of
    normal indices of ${\mathfrak s}$. Then
\begin{equation}\label{n0}
    {\mathfrak n}_1\le 2\kappa.
\end{equation}
 Moreover, if ${\mathfrak n}_1= 2\kappa$, then
\begin{equation}\label{sn0}
    {\mathfrak s}_{{\mathfrak n}_1-1}>0.
\end{equation}
\end{lem}
\begin{proof}
Since $\pi_\infty(F)\le\kappa$ it follows from~\eqref{infgzer}, that
\begin{equation}\label{GZNT}
     \infty \leq\lim_{\lambda\widehat{\rightarrow }\infty }
   {\lambda^{2\pi_{\infty}(F)+1}}F(\lambda) < 0.
\end{equation}
The normal index ${\mathfrak n}_1$ can be characterized by the
relations
\[
{\mathfrak s}_0=\dots={\mathfrak s}_{{\mathfrak
n}_1-2}=0,\,{\mathfrak s}_{{\mathfrak n}_1-1}\ne 0.
\]
Hence $F(\lambda)={\mathfrak F}(\lambda)/\lambda$ has the asymptotic expansion
\begin{equation}\label{asymp0}
{F}(\lambda)\sim -\frac{{\mathfrak s}_{{\mathfrak n}_1-1}}{\lambda^{{\mathfrak
n}_1+1}}-\dots-\frac{{\mathfrak s}_{2n}}{\lambda^{2n+2}}
-\dots,\quad\lambda\widehat{\rightarrow }\infty
\end{equation}
and \eqref{GZNT}  implies the inequality~\eqref{n0}.

If equality prevails in~\eqref{n0} then
$\pi_{\infty}(F)=\kappa$, the limit in \eqref{GZNT} is finite and
coincides with $-{\mathfrak s}_{{\mathfrak n}_1-1}$. This implies
the inequality~\eqref{sn0}.
\end{proof}
\begin{prop}
    Let ${\mathfrak F}\in {\bf D}_{\kk,-\infty}^\circ$ and $F(\lambda)={\mathfrak F}(\lambda)/{\lambda}$
    have asymptotic expansions~\eqref{asympF0} and~\eqref{eq:asymp}, and let
    $\cN_F$, $\cN({\bf s})$, $\cN({\mathfrak s})$ be defined by~\eqref{cN}, \eqref{NormInd1}, \eqref{NormIndFrak}. Then
    \[
\cN_F=\cN({\mathfrak s})\cap\cN({\bf s}).
    \]
\end{prop}
\begin{proof}
Let $\cN({\bf s})=\{n_j\}_{j=1}^\infty$. The statement is implied
by~\eqref{cN} and the equality
\[
\begin{vmatrix}
{\mathfrak s}_0 & \dots &  {\mathfrak s}_{n_{j}-1}\\
\hdotsfor{3}\\
{\mathfrak s}_{n_{j}-1}& \dots& {\mathfrak s}_{2n_{j}-2}
\end{vmatrix}
=\begin{vmatrix}
s_1 & \dots &  s_{n_{j}}\\
\hdotsfor{3}\\
s_{n_{j}}& \dots& s_{2n_{j}-1}
\end{vmatrix}\ne 0,
\]
which is immediate from~\eqref{sfrak} and~\eqref{cN2}.
\end{proof}

\subsection{The Schur transform of the ${\bf D}_{\kk,-\infty}^\circ$ functions.}
Let a  function ${\mathfrak F}\in {\bf D}_{\kk,-\infty}^\circ$ have
the asymptotic expansion~\eqref{asympF0}, let $\{{\mathfrak
n}_j\}_{j=1}^\infty$ be the set of normal indices for ${\mathfrak
s}=\{{\mathfrak s}_j\}_{j=0}^\infty$ and let ${\mathfrak
S}_n=({\mathfrak s}_{i+j})_{i,j=0}^n$. Let us set
 $   \epsilon_0=\mbox{sign }{\mathfrak s}_{{\mathfrak n}_1-1}$, 
\begin{equation}\label{p0}
    {\mathfrak p}_0(\lambda)=\frac{1 }{\det{\mathfrak S}_{{\mathfrak n}_1-1}}
    \det\left( \begin{array}{cccc}
    0           & \dots             & {\mathfrak s}_{{\mathfrak n}_1-1} & {\mathfrak s}_{{\mathfrak n}_1} \\
    \vdots      & \adots            & \adots                            & \vdots \\
    {\mathfrak s}_{{\mathfrak n}_1-1} & {\mathfrak s}_{{\mathfrak n}_1} & \dots                   & {\mathfrak s}_{{\mathfrak n}_1-1} \\
     1          & \lambda                & \dots                             & \lambda^{{\mathfrak n}_1} \\
                                                              \end{array}
                                                            \right)
    .
\end{equation}
The Schur transform of the function ${\mathfrak F}\in {\bf D}_{\kk,-\infty}^\circ$ is
defined by the equality
\begin{equation}\label{SchTr}
-\frac{1}{{\mathfrak    F}(\lambda)}=\epsilon_0{\mathfrak p}_0(\lambda)+    {\mathfrak
b}_0^2\wh{{\mathfrak F}}(\lambda),
\end{equation}
where ${\mathfrak b}_0$ is chosen in such a way that $\wh{{\mathfrak F}}$ has a
normalized expansion at $\infty$.
\begin{thm}\label{ScTr}
Let ${\mathfrak F}\in {\bf D}_{\kk,-\infty}^\circ$ and let $\wh{{\mathfrak F}}$ be the
Schur transform of ${{\mathfrak F}}$. Then:
\begin{enumerate}
  \item[(i)] $\wh{{\mathfrak F}}\in {\bf D}_{\kappa',-\infty}^\circ$ for some $\kappa'\le\kappa$;
  \item[(ii)] If ${{\mathfrak F}}\in {\bf D}_{1,-\infty}^\circ$ then $\wh{{\mathfrak F}}\in {\bf D}_{1,-\infty}^\circ$;
  \item[(iii)] The inverse Schur transform is given by
\begin{equation}\label{InvSchTr}
    {{\mathfrak F}}(\lambda)=-\frac{\epsilon_0}{{\mathfrak p}_0(\lambda)+\epsilon_0{\mathfrak b}_0^2\wh{{\mathfrak F}}(\lambda)}.
\end{equation}
\end{enumerate}
\end{thm}
\begin{proof}
(i) Direct calculations presented in~\cite[Lemmas 2.1, 2.4]{De} show
that $\wh{{\mathfrak F}}$ admits the asymptotic expansion
\begin{equation}\label{asympSTF}
\wh{{\mathfrak F}}(\lambda)\sim -\frac{{\mathfrak s}_{0}^{(1)}}{\lambda}-\frac{{\mathfrak
s}_{1}^{(1)}}{\lambda^{2}}-\dots-\frac{{\mathfrak s}_{2n}^{(1)}}{\lambda^{2n+1}}
-\dots,\quad\lambda\widehat{\rightarrow }\infty
\end{equation}
with some ${\mathfrak s}_{j}^{(1)}\in\dR$, $j\in\dZ_+$. Setting
\begin{equation}\label{eq:G1}
    G_1(\lambda):=\lambda\wh{{\mathfrak F}}(\lambda)
\end{equation}
one obtains from~\eqref{SchTr}
\begin{equation}\label{eq:G1Step}
    {\mathfrak b}_0^2G_1(\lambda)+\epsilon_0\lambda{\mathfrak p}_0(\lambda)=-\frac{1}{{
    F}(\lambda)}\in{\bf N}_{\kappa,-\infty}.
\end{equation}
Since $\deg \epsilon_0\lambda{\mathfrak p}_0(\lambda)={\mathfrak n}_{1}+1\ge
2$ then $\infty$ is a generalized pole of nonpositive type of the
polynomial $\epsilon_0\lambda{\mathfrak p}_0(\lambda)$ with multiplicity
\[
\kappa_\infty(\epsilon_0\lambda{\mathfrak p}_0(\lambda))\ge 1.
\]
It follows from~\eqref{asympSTF} and \eqref{eq:G1} that
\begin{equation}\label{G1s1}
    \lim_{\lambda\widehat{\rightarrow }\infty}G_1(\lambda)=-{\mathfrak s}_0^{(1)}
\end{equation}
and hence $\infty$ is not a generalized pole of nonpositive type of
$G_1$. By~\cite[Satz 1.13]{KL77} one obtains
\[
\kappa(G_1)+\kappa(\epsilon_0\lambda{\mathfrak
p}_0(\lambda))=\kappa(-1/F)=\kappa,
\]
and hence $G_1\in{\bf N}_{\kappa'',-\infty}$ for some
$\kappa''\le\kappa-1$.

Consider the function
\[
F_1(\lambda):=\frac{\wh{{\mathfrak F}}(\lambda)}{\lambda}=\frac{G_1(\lambda)}{\lambda^2}.
\]
It follows from~\eqref{infgzer} that the multiplicities of
generalized zeros at $\infty$ of $F_1$ and $G_1$ are related as
follows
\begin{equation}\label{eq:5.6}
\pi_\infty(F_1)=\pi_\infty(G_1)+1.
\end{equation}
So, by a theorem of M.G. Kre\u{\i}n and H. Langer~\cite[Theorem 3.5]{KL81}
$F_1\in{\bf N}_{\kappa',-\infty}$\,, where $\kappa'=\kappa''+1\le\kappa$.

(ii) By Proposition~\ref{lem:4.2} ${\mathfrak n}_1\le 2$ in the case
$\kappa=1$.

Assume first that ${\mathfrak n}_1=1$. Then $\deg \lambda{\mathfrak
p}_0(\lambda)=2$, $\kappa_\infty(\lambda{\mathfrak p}_0(\lambda))= 1$, and hence
$G_1\in{\bf N}_{0,-\infty}$. Then it follows from~\eqref{eq:5.6}
that $F_1\in{\bf N}_{1,-\infty}$.

Let now ${\mathfrak n}_1=2$. Then $\deg \lambda{\mathfrak p}_0(\lambda)=3$ and
in view of~\eqref{sn0} the leading coefficient of ${\mathfrak p}_0$
is positive. Therefore $\kappa_\infty(\lambda{\mathfrak p}_0(\lambda))= 1$,
and hence $G_1\in{\bf N}_{0,-\infty}$ and $F_1\in{\bf
N}_{1,-\infty}$.

(iii) The last statement is checked by straightforward calculations.
\end{proof}

\subsection{Diagonal Pad\'e approximants of the function ${\mathfrak F}\in
{\bf D}_{\kk,-\infty}$.} To prove the uniform convergence of diagonal
Pad\'e approximants for a function belonging to ${\bf D}_{\kk,-\infty}$, we need the
following lemma.
\begin{lem}[cf. \cite{DD07},\,\cite{Mag1}]\label{helplemma}
Let ${\mathfrak F}\in {\bf D}_{\kk,-\infty}$ and let
$F(\lambda):={\mathfrak F}(\lambda)/\lambda$. 
Then
\begin{equation}\label{shiftPA}
{{\mathfrak F}}^{[n-1/n-1]}(\lambda)=\lambda{F}^{[n/n-1]}(\lambda)\mbox{ for every }n\in\cN_F.
\end{equation}
\end{lem}
\begin{proof}
Suppose that $n\in\cN_F$. Then by Theorem~\ref{t_subdiagonal} the
Pad\'e approximant ${F}^{[n/n-1]}$ exists and
\begin{equation}\label{asforzero}
F^{[n/n-1]}(\lambda)+\sum_{j=0}^{2n-2}\frac{s_j}{\lambda^{j+1}}= O(\lambda^{-2n}),\quad \lambda\widehat{\rightarrow }\infty.
\end{equation}
Multiplying by $\lambda$ one obtains
\begin{equation}\label{shiftPA2}
\lambda F^{[n/n-1]}(\lambda)+\sum_{j=0}^{2n-2}\frac{s_j}{\lambda^{j}}= O(\lambda^{-(2n-1)}),
\quad \lambda\widehat{\rightarrow }\infty.
\end{equation}
Now the first term in~\eqref{shiftPA2} can be represented as
\[
\lambda F^{[n/n-1]}(\lambda)=\lambda\frac{A^{[n/n-1]}\left({1}/{\lambda}\right)}
{B^{[n/n-1]}\left({1}/{\lambda}\right)},
\]
where $\deg A^{[n/n-1]}\le n$, $\deg B^{[n/n-1]}\le n-1$, and $B^{[n/n-1]}(0)\ne 0$.
Moreover, it follows from the asymptotic expansion~\eqref{asforzero} that
$A^{[n/n-1]}(0)=0$. Hence, $A_1\left(\frac{1}{\lambda}\right)=\lambda
A^{[n/n-1]}\left(\frac{1}{\lambda}\right)$ is a polynomial in $\frac{1}{\lambda}$ of degree $\le
n-1$. This proves that
\[
\lambda F^{[n/n-1]}(\lambda)=\frac{A_1\left({1}/{\lambda}\right)} {B^{[n/n-1]}\left({1}/{\lambda}\right)},
\]
where $\deg A_1\le n-1$, $\deg B^{[n/n-1]}\le n-1$, and
$B^{[n/n-1]}(0)\ne 0$. So, it follows from~\eqref{shiftPA2} that
$\lambda F^{[n/n-1]}(\lambda)$ is the $[n-1/n-1]$ Pad\'e approximant for
${\mathfrak F}$.
\end{proof}
\begin{thm}\label{PAfordef}
Let ${\mathfrak F}\in {\bf D}_{\kk,-\infty}$ and let
\[
F(\lambda):=\frac{1}{\lambda}{\mathfrak F}(\lambda)\sim-\sum_{j=0}^\infty \frac{s_j}{{\lambda}^{j+1}},
\quad \lambda\widehat{\rightarrow }\infty,
\]
generate the determinate moment problem ${\bf M}_{\kk}({\bf s})$.
Then the sequence of diagonal Pad\'e approximants $\{{\mathfrak
F}^{[n-1/n-1]}\}_{n\in \cN_F}$ converges to ${\mathfrak F}$
locally uniformly on $\dC\setminus(\dR\cup\cP({\mathfrak F}))$.

Moreover, if the condition~\eqref{Pcond} is fulfilled for the
function $F$ of the form~\eqref{pertf} then the sequence of diagonal Pad\'e
approximants converges to ${\mathfrak F}$ locally uniformly
on $\dC\setminus
([-1-\varepsilon,1+\varepsilon]\cup\cP(\varphi))$ for some
$ \varepsilon>0$.
\end{thm}
\begin{proof}
It follows from Theorem~\ref{t_subdiagonal} that the sequence
$\{F^{[n/n-1]}\}_{n\in \cN_F}$
converges to $F$ locally uniformly on $\dC\setminus(\dR\cup\cP({
F}))$. Since $\dR\cup\cP({\mathfrak F})=\dR\cup\cP({F})$
the statement on the convergence on $\dC\setminus(\dR\cup\cP({\mathfrak F}))$
is implied by Lemma~\ref{helplemma}.

Under the condition~~\eqref{Pcond} the convergence on $\dC\setminus
([-1-\varepsilon,1+\varepsilon]\cup\cP(\varphi))$ for some
$ \varepsilon>0$ is a consequence of Theorem~\ref{th_subdA}.
\end{proof}
\begin{rem}
It should be noted that the above theorem and the appropriate variation of~\cite[Theorem~1.5.2]{Baker}
give us the possibility to make conclusions on
the locally uniform convergence of Pad\'e approximants for the function ${\mathfrak F}$
such that
\[
\frac{{\mathfrak F}(\lambda)}{\lambda+\zeta}\in{\bf
N}_{\kk,-\infty}
\]
for some $\zeta\in\dR$ and ${\mathfrak F}(\lambda)=O(1)$ as $\lambda\widehat{\rightarrow }\infty$.
\end{rem}

\subsection{Generalized Jacobi matrix associated with the function ${\mathfrak F}\in
{\bf D}_{\kk,-\infty}^\circ$.}
\begin{thm}\label{GJM}
Let ${\mathfrak F}\in {\bf D}_{\kk,-\infty}^\circ$,  let the
sequence ${\mathfrak s}=\{{\mathfrak s}_{j}\}_{j=0}^\infty$ be
defined by the asymptotic expansion~\eqref{asympF0},
    and let $ \cN({\mathfrak s})=\{{\mathfrak n}_i\}_{i=1}^\infty\,$ be the set of
    normal indices of ${\mathfrak s}$.
Then:
\begin{enumerate}
  \item[(i)] ${\mathfrak F}$ admits the expansion into the $P$-fraction
\begin{equation}\label{PfractionG}
-\frac{\epsilon_0}{{\mathfrak p}_0(\lambda)}
\begin{array}{l} \\ -\end{array}
\frac{\epsilon_0\epsilon_1{\mathfrak b}_0^2}{{\mathfrak p}_1(\lambda)}
\begin{array}{ccc} \\ - & \cdots & -\end{array}
\frac{\epsilon_{N-1}\epsilon_N {\mathfrak b}_{N-1}^2}{{\mathfrak
p}_N(\lambda)}\begin{array}{l}\\-\dots\end{array},
\end{equation}
where
  ${\mathfrak p}_i$ are polynomials of degree ${\mathfrak k}_i:={\mathfrak n}_{i+1}-{\mathfrak n}_{i}$
  $(\le 2\kappa)$, $\epsilon_i=\pm 1
  $, ${\mathfrak b}_i>0$, $i\in\dZ_+$;
  \item[(ii)] If $\sJ$ is the generalized Jacobi matrix associated with the $P$-fraction~\eqref{PfractionG},
  and ${\mathfrak P}_i(\lambda)$, ${\mathfrak Q}_i(\lambda)$ are given by
  \begin{eqnarray}
{\mathfrak P}_{i}(\lambda)&=&({\mathfrak b}_{0}\dots{\mathfrak b}_{i-1})^{-1}\det(\lambda-{\mathfrak J}_{[0,i-1]}),\label{Wtpolynom1}\\
{\mathfrak Q} _{i}(\lambda)&=&({\mathfrak b}_{0}\dots {\mathfrak
b}_{i-1})^{-1}\det(\lambda-{\mathfrak J}_{[1,i-1]}),\label{Wtpolynom2}
\end{eqnarray}
then the $i$-th convergent to~\eqref{PfractionG} coincides with
$-{\mathfrak Q}_i(\lambda)/{\mathfrak P}_i(\lambda)$ and is the
$[{\mathfrak n}_i/{\mathfrak n}_i]$ Pad\'e approximant for
${\mathfrak F}(\lambda)$.
\end{enumerate}
\end{thm}
\begin{proof}
(i) It follows from Theorem~\ref{ScTr} that any function ${\mathfrak
F}\in {\bf D}_{\kk,-\infty}^\circ$ can be represented as follows
\[
{{\mathfrak F}}(\lambda)=-\frac{\epsilon_0}{{\mathfrak
p}_0(\lambda)+\epsilon_0{\mathfrak b}_0^2{{\mathfrak F}_1}(\lambda)},
\]
where ${\mathfrak p}_0$ is a monic polynomial of degree $\deg
{\mathfrak p}_0={\mathfrak n}_1\le 2\kappa$ (see
formula~\eqref{n0}), $\epsilon_0=\pm 1$, $b_0\in\dR_+$, and
${\mathfrak F}_1\in{\bf D}_{\kk_1,-\infty}^\circ$ with
$\kk_1\le\kk$. Further, one can apply Theorem~\ref{ScTr} to
${\mathfrak F}_1$ and so on. Thus, the Schur algorithm leads
to~\eqref{PfractionG}. To complete the proof, note that the relation
$\deg {\mathfrak p}_i={\mathfrak n}_{i+1}-{\mathfrak n}_{i}$ follows
from~\cite[Corollary 3.6]{De}.

\noindent (ii) This part is proved in line
with~\cite[Proposition 2.3]{DD} (see also~\cite{DD07}).
\end{proof}

\begin{rem}
Let $\wh{{\mathfrak P}}_{i}$ 
be monic polynomials
associated with  ${\mathfrak P}_{i}$ 
by the equalities
\[
\wh{{\mathfrak P}}_{i}(\lambda):=({\mathfrak b}_{0}\dots{\mathfrak
b}_{i-1}){\mathfrak
P}_{i}(\lambda),\quad \mbox{deg }{\mathfrak P}_{i}={\mathfrak n}_{i},\quad i\in\dN.
\]
Then it follows from Theorem~\ref{th_subdA}, Lemma~\ref{helplemma}, and
Theorem~\ref{GJM} that $\wh{{\mathfrak P}}_{i}$ 
are the Christoffel transformations of the polynomials $\wh P_j$
corresponding to $F(\lambda)={\mathfrak F}(\lambda)/\lambda$
\[
\begin{split}
\wh{{\mathfrak P}}_{i}(\lambda)&
=(b_0\dots
b_{j-1})\left(P_j(\lambda)-\frac{P_j(0)}{P_{j-1}(0)}P_{j-1}(\lambda)\right)\frac{1}{\lambda}\\
&=\frac{\wh P_j(\lambda)\wh P_{j-1}(0)-\wh P_{j-1}(\lambda)\wh P_{j}(0)}{\wh P_{j-1}(0)\lambda},\\
\end{split}
\]
such that ${{\mathfrak n}}_{i}=\deg \wh{{\mathfrak P}}_{i}=\deg \wh{
P}_{j}=n_j-1$, $n_j\in\cN({\mathbf s})$. Since $F\in{\mathbf
N}_\kappa$, the Christoffel transformation is defined for every
natural $n(=n_j)$ large enough.
\end{rem}

In the case when ${\mathfrak F}\in {\bf D}_{1,-\infty}^\circ$ one
can simplify the form of the generalized Jacobi matrix ${\mathfrak
J}$.
\begin{prop}\label{GJM0}
Let ${\mathfrak F}\in {\bf D}_{1,-\infty}^\circ$ satisfy the
assumptions of Theorem~\ref{GJM} and let ${\mathfrak J}$ be the
generalized Jacobi matrix associated with the
$P$-fraction~\eqref{PfractionG}. Then ${\mathfrak k}_i:={\mathfrak
n}_{i+1}-{\mathfrak n}_{i}$ is either 1 or 2 and the block matrix
${\mathfrak A}_i$ in ${\mathfrak J}$ takes the form
\[
{\mathfrak A}_i=\left\{\begin{array}{cc}
             -{\mathfrak p}_0^{(i)}, & \mbox{if }{\mathfrak k}_i=1; \\
             \begin{pmatrix} 0 & -{\mathfrak p}_0^{(i)}\\
              1 & -{\mathfrak p}_1^{(i)}\end{pmatrix} & \mbox{if }{\mathfrak k}_i=2,
           \end{array}\right.
\]
where ${\mathfrak p}_0^{(i)}$, ${\mathfrak p}_1^{(i)}$ are
coefficients of the polynomials ${\mathfrak p}_i$
in~\eqref{PfractionG}.
\end{prop}
It may happen that the generalized Jacobi matrix ${\mathfrak J}$ is
unbounded even in the case when the support of ${\mathfrak F}$ is
bounded (see Example~\ref{Ex:1}). It should be noted that bounded
generalized Jacobi matrices associated to~\eqref{PfractionG} were
considered in~\cite{De09}.

\section{Particular cases}
\label{sect5}
\subsection{The case when ${\mathfrak F}$ is holomorphic at $\infty$.}
Consider the function ${\mathfrak F}$ of the form
\begin{equation}\label{pertf1}
{\mathfrak F}(\lambda)=r_1(\lambda)\int_{-1}^1\frac{td\sigma(t)}{t-\lambda}+r_2(\lambda),
\end{equation}
where $\sigma$, $r_1$  and $r_2$ satisfy the assumptions (A1), (A2)
and (A$3'$). 
\begin{thm}\label{theoremA}
Let the function ${\mathfrak F}$ be of the form~\eqref{pertf1} and
let $F(\lambda)={\mathfrak F}(\lambda)/\lambda$. Then ${\mathfrak
F}\in {\bf D}_{\kappa,-\infty}$ for some $\kappa\in\dZ_+$ and the
sequence of $[n-1/n-1]$ Pad\'e approximants $\{{\mathfrak
F}^{[n-1/n-1]}\}_{n\in\cN_F}$ converges to ${\mathfrak F}$ locally
uniformly in $\dC\setminus (\dR\cup\cP(\varphi))$.

Moreover, if the condition~\eqref{Pcond} is fulfilled for the
function $F$ then the sequence $\{{\mathfrak
F}^{[n-1/n-1]}\}_{n\in\cN_F}$
converges to ${\mathfrak F}$ locally uniformly in $\dC\setminus
([-1-\varepsilon,1+\varepsilon]\cup\cP(\varphi))$ for some $
\varepsilon>0$.
\end{thm}
\begin{proof}
The function ${\mathfrak F} (\lambda)/\lambda$ admits the representation
\[ F(\lambda)=\frac{{\mathfrak
F}(\lambda)}{\lambda}=r_1(\lambda)\int_a^b\frac{d\sigma(t)}{t-\lambda}+\wt r_2(\lambda),
\]
where
\[
\wt r_2(\lambda)=\frac{r_1(\lambda)}{\lambda}\int_a^b{d\sigma(t)}+\frac{r_2(\lambda)}{\lambda}.
\]
Therefore $F\in{\bf N}_{\kappa,-\infty}$ (see~\cite{KL77}), and hence ${\mathfrak F}\in
{\bf D}_{\kappa,-\infty}$.

The statements concerning convergence of the sequence of
diagonal Pad\'e approximants of ${\mathfrak F}$ are implied by
Theorem~\ref{PAfordef}.
\end{proof}
\begin{rem}\label{Dkappa}
    In fact, as is easily seen from~\cite{KL77}, every function ${\mathfrak F}\in {\bf
D}_{\kappa,-\infty}$ holomorphic at infinity admits the representation~\eqref{pertf1}.
\end{rem}

\begin{ex}\label{Ex:1}
Let $\theta\in\dR$ be an irrational number and consider the function
\[
{\mathfrak
F}(\lambda)=\int_{-1+\cos{\pi\theta}}^{1+\cos{\pi\theta}}\frac{td\sigma(t)}{t-\lambda},
\mbox{where }d\sigma(t)=\frac{dt}{\sqrt{1-(t-\cos{\pi\theta})^2}}.
\]
Substitution $x=t-\cos{\pi\theta}$ leads to the equality
\[
{\mathfrak
F}(\lambda+\cos{\pi\theta})=\int_{-1}^{1}\frac{(x+\cos{\pi\theta})d\omega(x)}{x-\lambda},
\mbox{where }d\omega(x)=\frac{dx}{\sqrt{1-x^2}}.
\]
As was shown in~\cite{Stahl} every point of $\dR$ is an accumulation
point of the set of poles of the diagonal Pad\'e approximants for ${\mathfrak
F}(\cdot+\cos{\pi\theta})$. As a consequence, the diagonal
Pad\'e approximants for ${\mathfrak F}$ do not converge on
$\dR\setminus[-1+\cos{\pi\theta},1+\cos{\pi\theta}]$. Therefore, the
corresponding generalized Jacobi matrix ${\mathfrak J}$ is
unbounded.

However,  there exists a subsequence of ${\mathfrak F}^{[n-1/n-1]}$
converging  in a neighborhood of $\infty$. Indeed, applying
Lemma~\ref{helplemma} to the function ${\mathfrak F}-\gamma$ with
$\gamma=\int\limits_{1+\cos{\pi\theta}}^{-1+\cos{\pi\theta}}d\sigma(t),$
one obtains
\[
{\mathfrak
F}^{[n-1/n-1]}(\lambda)=\lambda{F}^{[n/n-1]}(\lambda)+\gamma, \quad
\mbox{where }
F(\lambda)=\int\limits_{-1+\cos{\pi\theta}}^{1+\cos{\pi\theta}}\frac{d\sigma(t)}{t-\lambda}.
\]
Clearly, the shifted Chebyshev polynomials
$T_n(\cdot-\cos{\pi\theta})$ are orthonormal  with respect to
$\sigma$. Consequently, we can calculate explicitly the coefficient
\[
\tau_n=-\frac{T_{n+1}(\cos{\pi\theta})}{2T_n(\cos{\pi\theta})}=-
\frac{\cos(n+1)\pi\{|\theta|\}}{2\cos n\pi\{|\theta|\}}=\frac{1}{2}
(\cos\pi\{|\theta|\}-\sin\pi\{|\theta|\}\tan n\pi\{|\theta|\}),
\]
where $n\in\dN$ and $\{x\}$ denotes the fractional part of
$x\in\dR$. Since the set $\{\{n|\theta|\}\}_{n=0}^{\infty}$ is dense
in $(0,1)$,  there is a bounded subsequence of
$\{\tau_n\}_{n=0}^{\infty}$ and thus, by Remark~\ref{Subs},  there
exists a subsequence of diagonal Pad\'e approximants converging in a
neighborhood of $\infty$.
\end{ex}

\begin{rem}
Let us consider a function ${\mathfrak F}$ of the following form
\begin{equation}\label{fmodPA}
{\mathfrak F}(\lambda)=\int_{-1}^{1}\frac{td\sigma(t)}{t-\lambda},
\end{equation}
where $\sigma$ is a nonnegative probability measure on $[-1,1]$.
It is clear that
\[
{\mathfrak F}(\lambda)=(\cJ(\cJ-\lambda)^{-1}e,e)_{\ell^2_{[0,\infty)}}
=1+\lambda((\cJ-\lambda)^{-1}e,e)_{\ell^2_{[0,\infty)}},
\]
where $\cJ$ is the classical Jacobi matrix constructed by the measure $\sigma$ via
the usual procedure~\cite{Ach61}.
Now, let us consider the following modified Pad\'e approximant
\begin{equation}\label{modPA}
{\mathfrak F}_{*}^{[n/n]}(\lambda)=(\cJ_{[0,n-1]}(\cJ_{[0,n-1]}-\lambda)^{-1}e,e)=
\frac{P_n(\lambda)-\lambda Q_n(\lambda)}{P_n(\lambda)}
\end{equation}
where $P_n$, $Q_n$ are polynomials of the first and second kinds
corresponding to the measure $\sigma$. It follows from the
Markov theorem (as well as from the spectral decomposition theorem)
that
\[
{\mathfrak F}_{*}^{[n/n]}\to {\mathfrak F}
\]
locally uniformly in $\dC\setminus[-1,1]$. So,  to avoid the
phenomenon described in the above example one can use the modified
Pad\'e approximants~\eqref{modPA} for the function ${\mathfrak F}$
of the form~\eqref{fmodPA}.
\end{rem}
\subsection{The case when $\mbox{supp }\sigma$ has a gap.}
Assume now that $r_1(\lambda)\equiv r_2(\lambda)\equiv 1$ in~\eqref{pertf1} and the support $E$ of
the finite nonnegative measure $\sigma$  is contained in the union of two intervals
\[
E=[-1,\alpha]\cup[\beta,1],\quad \alpha<0<\beta.
\]
First, we will show that in this case the diagonal Pad\'e approximants for
${\mathfrak F}$ have no poles inside the gap
$(\alpha,\beta)$.
\begin{prop}\label{GJM2}
Let $\sigma$ be a finite nonnegative measure on $E=[-1,\alpha]\cup[\beta,1]$ and let
\begin{equation}\label{PartCase2}
    {\mathfrak F}(\lambda)=\int_E\frac{td\sigma(t)}{t-\lambda}.
\end{equation}
Then:
\begin{enumerate}
\item[(i)] ${\mathfrak F}\in {\bf
D}_{1,-\infty}^\circ$;
\item[(ii)]
The polynomials  ${\mathfrak P}_j$, $j\in\dN$ have no zeros inside the gap
$(\alpha,\beta)$;
\item[(iii)]
The function
\begin{equation}\label{PartCase20}
    {\mathfrak F}_0(\lambda)=\int_E\frac{td\sigma(t)}{t-\lambda}-\gamma,\quad \gamma=\int_E{d\sigma(t)}.
\end{equation}
belongs to the class ${\bf D}_{0,-\infty}$.
\end{enumerate}
\end{prop}
\begin{proof}
(i) The first statement is implied by the equality
\begin{equation}\label{psi}
F(\lambda)=\frac{{\mathfrak
F}(\lambda)}{\lambda}=\int_E\frac{d\sigma(t)}{t-\lambda}+\frac{1}{\lambda}\int_E d\sigma(t),
\end{equation}
since $F\in {\bf N}_{1,-\infty}$.

(ii) Next, ${\mathfrak P}_j(0)$ coincides with the Hankel determinants
\[
{\mathfrak P}_j(0)=\left|\begin{array}{ccc}
               {\mathfrak s}_1 & \dots & {\mathfrak s}_j\\
               \hdotsfor{3} \\
               {\mathfrak s}_j & \dots & {\mathfrak s}_{2j-1}
             \end{array}
\right|
\]
which are positive since $s_j={\mathfrak s}_{j+1}$, $ j=0,1,\dots$ are moments of the
positive measure $t^2d\sigma(t)$ with infinite support $E$.

Similarly, ${\mathfrak P}_j(\theta)$ coincides with the Hankel determinants
\[
{\mathfrak P}_j(\theta)=\left|\begin{array}{ccc}
               {\mathfrak s}_1-\theta{\mathfrak s}_0 & \dots & {\mathfrak s}_j-\theta{\mathfrak s}_{j-1}\\
               \hdotsfor{3} \\
               {\mathfrak s}_j-\theta{\mathfrak s}_{j-1} & \dots & {\mathfrak s}_{2j-1}-\theta{\mathfrak s}_{2j-2}
             \end{array}
\right|
\]
which are positive for $\theta\in(\alpha,\beta)$, since $s_j:={\mathfrak
s}_{j+1}-\theta{\mathfrak s}_{j}$ are moments of the positive measure
$t(t-\theta)d\sigma(t)$.

(iii) The third statement follows from the equality
\begin{equation}\label{psi1}
F_0(\lambda)=\frac{{\mathfrak F}_0(\lambda)}{\lambda}=\int_E\frac{d\sigma(t)}{t-\lambda}.
\end{equation}
\end{proof}

Next, one can apply Theorem~\ref{th_subdA} to prove the convergence of
diagonal Pad\'e approximants for ${\mathfrak F}$ on the real line.
\begin{thm}\label{theoremC}
Let $\sigma$ be a finite nonnegative measure on $E=[-1,\alpha]\cup[\beta,1]$,  let
${\mathfrak F}$ have the form~\eqref{PartCase2} and let $\{\wh P_j\}_{j=0}^{\infty}$
be  the set of normalized polynomials orthogonal with respect to $\sigma$. Then:
\begin{enumerate}
\item[(i)] The sequence of diagonal Pad\'e approximants
$\{{\mathfrak F}^{[n-1/n-1]}\}_{n\in\cN_F}$
converges to ${\mathfrak F}$ locally
uniformly in $\dC\setminus((-\infty,\alpha]\cup[\beta,\infty))$;

\item[(ii)] The sequence $\{{\mathfrak
F}^{[n-1/n-1]}\}_{n\in\cN_F}$ 
converges to ${\mathfrak F}$ locally uniformly in \linebreak
$\dC\setminus ([-1-\varepsilon,\alpha]\cup[\beta, 1+\varepsilon])$
for some $ \varepsilon>0$ if and only if the
condition~\eqref{WhPcond} is fulfilled.
\end{enumerate}
\end{thm}
\begin{proof}
(i) As follows from Lemma~\ref{helplemma} and~\cite[Theorem
4.16]{DD07} the Pad\'e approximant ${\mathfrak F_0}^{[n-1/n-1]}$
takes the form
\begin{equation}\label{eq:5.4}
{\mathfrak F_0}^{[n-1/n-1]}(\lambda)=\lambda {
F}_0^{[n/n-1]}(\lambda)=\lambda
((\cJ_{[0,n-1]}^{(K)}-\lambda)^{-1}e,e),\quad n=n_j\in\cN_{F_0},
\end{equation}
where $\cJ$ is a classical Jacobi matrix corresponding to the
measure $\sigma$, and $\cJ_{[0,n-1]}^{(K)}$ is defined in
Remark~\ref{zero_jacobi2}. Let us emphasize that
$\cJ_{[0,n-1]}^{(K)}$ is a classical Jacobi matrix since
$F_0\in{\mathbf N}_0$. By Proposition~\ref{GJM2} and
Theorem~\ref{GJM} ${\mathfrak F_0}^{[n-1/n-1]}$ is holomorphic on
$(\alpha,0)\cup(0,\beta)$ and hence by~\eqref{eq:5.4} and
Remark~\ref{ResSet} the set $(\alpha,0)\cup(0,\beta)$ is contained
in $\rho(\cJ_{[0,n-1]})$. It follows from the spectral theorem that
for arbitrary $\varepsilon>0$
\[
\|(\cJ_{[0,n-1]}^{(K)}-\lambda)^{-1}\|\le\frac{1}{\varepsilon}\quad\mbox{for all
}\lambda\in(\alpha+\varepsilon,-\varepsilon)\cup(\varepsilon,\beta-\varepsilon).
\]
Then by the Vitali theorem the sequence
$\{{\mathfrak F_0}^{[n-1/n-1]}(\lambda)\}_{n\in\cN_F}$
converges to
\[
{\mathfrak F_0}(\lambda)=\lambda ((\cJ-\lambda)^{-1}e,e)
\]
locally uniformly on $\dC\setminus((-\infty,\alpha]\cup\{0\}\cup[\beta,\infty))$.
Moreover, for $\varepsilon>0$ small enough ${\mathfrak F_0}^{[n-1/n-1]}(\lambda)$
converges to ${\mathfrak F_0}(\lambda)$ uniformly on the circle $|\lambda|=\varepsilon$.
Then by the mean value theorem
\[
{\mathfrak F_0}^{[n-1/n-1]}(0)=\frac{1}{2\pi}\int_0^{2\pi}{\mathfrak
F_0}^{[n-1/n-1]}(\varepsilon e^{it})dt\to \frac{1}{2\pi}\int_0^{2\pi}{\mathfrak
F_0}(\varepsilon e^{it})dt={\mathfrak F_0}(0)
\]
as $n\to\infty$. To complete the proof it remains to mention that
the $[n-1/n-1]$ Pad\'e approximants  ${\mathfrak F}^{[n-1/n-1]}$,
$n\in\cN_F$ are connected with the $[n-1/n-1]$ Pad\'e approximants
${\mathfrak F}_0^{[n-1/n-1]}$ by the equality
\[
{\mathfrak F}^{[n-1/n-1]}={\mathfrak F}_0^{[n-1/n-1]}+\gamma.
\]

(ii) The necessity of the second statement is contained in
Theorem~\ref{PAfordef}. Let us prove the sufficiency by proving
the inverse statement. So, suppose that
\begin{equation}
\sup\limits_{n_j\in \cN_F}\left|\frac{\wh P_{j}(0)}{\wh
P_{j-1}(0)}\right|=\infty.
\end{equation}
Due to~\eqref{eq:5.4} the poles of the Pad\'e approximant
${\mathfrak F_0}^{[n-1/n-1]}$ coincide with eigenvalues of the
matrix $\cJ_{[0,n-1]}^{(K)}$. Taking into account that
$\cJ_{[0,n-1]}^{(K)}$ is a self-adjoint matrix, one obtains
\[
|\lambda_{max}(\cJ_{[0,n-1]}^{(K)})|=\Vert \cJ_{[0,n-1]}^{(K)}\Vert
\ge|(\cJ_{[0,n-1]}^{(K)}e_{n-1},e_{n-1})|=\left|p_0^{(n-1)}-\frac{\wh
P_{n}(0)}{\wh P_{n-1}(0)}\right|,
\]
where $\lambda_{max}(\cJ_{[0,n-1]}^{(K)})$ is the eigenvalue of the
matrix $\cJ_{[0,n-1]}^{(K)}$ with the largest absolute value.
Since the sequence $\{p_0^{(n-1)}\}_{n=1}^{\infty}$ is bounded, we
have that infinity is an accumulation point of the set of all
poles of Pad\'e approximants ${\mathfrak F_0}^{[n-1/n-1]}$.
\end{proof}

\begin{rem}
Let ${\mathfrak F}$ be a function having the following form
\[
{\mathfrak F}(\lambda)=\int_{-1}^{1}\frac{\rho(t)dt}{t-\lambda},
\]
where $\rho$ is a nonvanishing on $[-1,1]$ complex-valued function.
Under some assumptions on $\rho$, the locally uniform convergence of
the diagonal Pad\'e approximants for ${\mathfrak F}$ was proved by
A.~Magnus~\cite{Mag87} (see also~\cite{Stahl}). Using the technique
of Riemann-Hilbert problems, the result was reproved by
A.I.~Aptekarev and W.~van~Assche~\cite{Apt}.
\end{rem}

\begin{ex}
Let ${\mathfrak F}$ have the form~\eqref{PartCase2} with
an absolutely continuous measure $d\sigma(t)=\rho(t)dt$ such that
$\rho(t)$ is an even function on $E=[-1,-\beta]\cup[\beta,1]$ and
$\rho(t)=0$ for $t\in\dR\setminus E$. Then the polynomials $P_{2j+1}$
are odd (for instance, see~\cite[formula~(5.90)]{Si})
and, therefore, $P_{2j+1}(0)=0$ for all $j\in \dN$. Hence the
condition~\eqref{Pcond} is fulfilled and by Theorem~\ref{theoremC}
the Pad\'e approximants ${\mathfrak F}^{[2j/2j]}$ converge to
${\mathfrak F}$ on $\dC\setminus E$. This fact can be shown
directly, since ${\mathfrak F}(\lambda)$ admits the representation
\[
{\mathfrak F}(\lambda)=\varphi(\lambda^2),\quad \mbox{where }
\varphi(\mu)=\int_{\beta^2}^1\frac{\sqrt{s}\rho(\sqrt{s})ds}{s-\mu}\in
{\bf N}_0
\]
and, hence,
 $ {\mathfrak
F}^{[2j/2j]}(\lambda)=\varphi^{[j/j]}(\lambda^2) $ converge to ${\mathfrak
F}(\lambda)=\varphi(\lambda^2)$ for all $\lambda^2\in\dC\setminus[\beta^2,1]$,
or, equivalently, for all $\lambda\in\dC\setminus E$.
\end{ex}

Due to Remark~\ref{Subs} it is enough to find a bounded subsequence
of $\{\tau_{j-1}\}_{n_j\in\cN_F}$ to say that there exists a
subsequence of diagonal Pad\'e approximants of ${\mathfrak F}$ which
converges locally uniformly in a neighborhood of $\infty$. In the
following proposition we find a sufficient condition for the
boundedness of a subsequence of $\{\tau_{j-1}\}_{n_j\in\cN_F}$.

\begin{prop}\label{P_ac}
Let $\sigma$ be a finite nonnegative measure on $E=[-1,\alpha]\cup[\beta,1]$  and let
$\{\wh P_j\}_{j=0}^{\infty}$ be  the set of monic polynomials orthogonal with respect to
$\sigma$. Assume that 0 is not an accumulation point of zeros of a subsequence $\{\wh
P_{j_k}\}_{k=1}^{\infty}$. Then the sequence $\left\{\frac{\wh P_{j_k+1}(0)}{\wh
P_{j_k(0)}}\right\}_{k=1}^{\infty}$ is bounded.
\end{prop}
\begin{proof}
The orthogonal polynomials $\wh P_j$ satisfy the following recurrence relations
\[
\lambda \wh P_j(\lambda)=b_{j-1}^2\wh P_{j-1}(\lambda)+a_j\wh P_{j}(\lambda)+\wh
P_{j+1}(\lambda),\quad b_j>0,\quad a_j\in\dR,
\]
which implies the following equality
\begin{equation}\label{Mfun}
\frac{\wh P_{j+1}(\lambda)}{\wh P_j(\lambda)}=\lambda-a_j+b_{j-1}^2\left(-\frac{\wh
P_{j-1}(\lambda)}{\wh P_j(\lambda)}\right).
\end{equation}
It is well known (see~\cite{Akh}) that $-\frac{\wh P_{j-1}(\lambda)}{\wh P_j(\lambda)}$
belongs to ${\bf N}_0$ and hence there is a nonnegative measure $\sigma^{(j)}$ such
that
\[
-\frac{\wh P_{j-1}(\lambda)}{\wh
P_j(\lambda)}=\int_{-1}^{1}\frac{d\sigma^{(j)}(t)}{t-\lambda}.
\]
Moreover $\sigma^{(j)}$ satisfies the condition
\[\int_{-1}^{1}d\sigma^{(j)}(t)=1\]
because of the asymptotic relation $-\frac{\wh P_{j-1}(\lambda)}{\wh P_j(\lambda)}=
-\frac{1}{\lambda}+o\left(\frac{1}{\lambda}\right)$ as $\lambda\to\infty$. Since the zeros of
$\{\wh P_{j_k}\}_{k=1}^{\infty}$ do not
accumulate to 0, there exists $\delta>0$ such that
\[
(-\delta,\delta)\cap\supp \sigma^{(j_k)}=\emptyset,\quad
k=1,2,\dots
\]
So, we have the following estimate
\begin{equation}\label{Mf_est}
\left|\frac{\wh P_{j-1}(0)}{\wh
P_j(0)}\right|=\left|\int_{-1}^{1}\frac{d\sigma^{(j_k)}(t)}{t}\right|\le\frac{1}{\delta}.
\end{equation}
Finally, the boundedness of $\left\{\frac{\wh P_{j_k+1}(0)}{\wh
P_{j_k(0)}}\right\}_{k=1}^{\infty}$ follows from the boundedness of the sequences
$\{a_j\}_{j=0}^{\infty}$, $\{b_j\}_{j=0}^{\infty}$, the estimate~\eqref{Mf_est}, and the
equality~\eqref{Mfun}.
\end{proof}

\subsection{The case when $\sigma$ satisfies the Szeg\"o conditions.}
Now, a natural question arises: under what conditions 0 is not an accumulation point of
zeros of a subsequence $\{P_{j_k}\}_{j=0}^{\infty}$? In this subsection the answer to
this question is given for functions of the form~\eqref{PartCase2} under an additional
assumption that the measure $d\sigma(t)=\rho(t)dt$ satisfies the Szeg\"o condition on
each of the intervals $[-1,\alpha]$ and $[\beta,1]$
\begin{equation}\label{Szego}
    \int_{-1}^{\alpha}\frac{\log\rho(t)}{\sqrt{(\alpha-t)(t+1)}}dt>-\infty,\quad
    \int_{\beta}^{1}\frac{\log\rho(t)}{\sqrt{(1-t)(t-\beta)}}dt>-\infty.
\end{equation}

As is known the polynomials $P_j$ have at most one zero in the
interval $(\alpha,\beta)$. The information about accumulation
points of these zeros can be formulated in terms of the harmonic
measure $\omega(\lambda)$ of $[-1,\alpha]$ with respect to
$\dC\setminus E$, i.e. harmonic function on $\dC\setminus E$ whose
boundary values are equal 1 on $[-1,\alpha]$ and 0 on $[\beta,1]$.
\begin{rem}
For more detailed and deep analysis of the behavior of zeros of
orthogonal polynomials see~\cite{Suet2002} (see
also~\cite{Suet2005}).
\end{rem}

Assume first that $\omega(\infty)$ is an irrational number. Then
by a theorem of E.~Rakh\-manov (\cite[Theorem 0.2]{Rak}) every
point of $(\alpha,\beta)$ and, in particular, 0 is an accumulation
point of zeros of a sequence $\{P_{j}\}_{j=0}^{\infty}$. However, since there is
only one zero of $P_j$ in the gap $(\alpha,\beta)$ it is possible to choose
a subsequence of $\{P_{j}\}_{j=0}^{\infty}$ which zeros do not accumulate to 0.
Further, as follows from Proposition~\ref{P_ac} there is a
subsequence of $[n/n]$ Pad\'e approximants of ${\mathfrak F}$
which converges to ${\mathfrak F}$ locally uniformly on
$\dC\setminus([-1-\varepsilon,\alpha]\cup[\beta,1+\varepsilon])$
for some $\varepsilon>0$.

Assume now that $\omega(\infty)$ is a rational number $m/n$, where $m,n\in\dN$ and
$gcd(m,n)=1$. Then it follows from~\cite[formula (57)]{Rak} that  every accumulation point of
zeros of polynomials $\{P_{j}\}_{j=1}^{\infty}$ in the interval $(\alpha,\beta)$
satisfies one of the equation
\begin{equation}\label{R_eq}
\omega_1(z)\equiv \frac{k}{n}(\mbox{mod }2), \quad k\in\dZ, |k|\le n.
\end{equation}
In the case when 0 is not a solution of the equation~\eqref{R_eq}
it follows from Proposition~\ref{P_ac} and Theorem~\ref{theoremC}
that the sequence $\{{\mathfrak F}^{[n_j-1/n_j-1]}\}_{n_j\in\cN_f}$
of  Pad\'e approximants of ${\mathfrak F}$ converges to
${\mathfrak F}$ locally uniformly on
$\dC\setminus([-1-\varepsilon,\alpha]\cup[\beta,1+\varepsilon])$
for some $\varepsilon>0$.

The harmonic measure of $[-1,\alpha]$ with respect to $\dC\setminus E$ can be calculated
explicitly (see~\cite{Akh}). Let real $k$ be defined by
\[
k^2=\frac{2(\beta-\alpha)}{(1-\alpha)(1+\beta)}.
\]
Consider the function $x=\mbox{sn }w$ with the modulus $k$, and with the primitive
periods $4K$, $2iK'$. As is known (see~\cite[p.190]{Akh}) the mapping
\begin{equation}\label{eq:6.3}
    z=\alpha+\frac{1-\alpha^2}{2\mbox{sn}^2\frac{K'\ln w}{\pi}+\alpha-1}
\end{equation}
maps conformally the ring
\begin{equation}\label{Ring}
r:=e^{-\frac{\pi K}{K'}}<|w|<1
\end{equation}
onto the plane $\dC$ with the cuts $[-1,\alpha]\cup[\beta,1]$, moreover, the semicircle
$|w|=1$ $(\mbox{Im }w\ge 0)$ is mapped onto the upper shore of the cut $[-1,\alpha]$.

As is well known (see~\cite{Land})  the harmonic measure $\omega_R$ of the ring
\eqref{Ring} has  the form
\[
\omega_R(w)=\frac{\ln |w|-\ln r}{\ln 1-\ln r}=\frac{\ln |w|+\frac{\pi K}{K'}}{\frac{\pi
K}{K'}}.
\]
So, the harmonic measure of $[-1,\alpha]$ with respect to $\dC\setminus E$ can be found by
\begin{equation}\label{eq:6.4}
\omega(z)=\frac{K'}{\pi K}\ln |w|+1.
\end{equation}

Let us choose $w_\infty\in(r,1)$ and $w_0\in(-1,r)$ such that
\begin{equation}\label{eq:6.5}
1-2\mbox{sn}^2\frac{K'\ln w_\infty}{\pi}=\alpha,\quad 1-2\mbox{sn}^2\frac{K'\ln
w_0}{\pi}=\frac{1}{\alpha}.
\end{equation}
Then the  numbers $w_\infty$,  $w_0$ correspond to  $z=\infty$ and $z=0$
via~\eqref{eq:6.3}. It follows from~\eqref{eq:6.3} that $\omega(\infty)$ is a rational
number $\frac{m}{n}$ if and only if $\frac{K'}{\pi K}\ln w_\infty=\frac{m-n}{n}$, which
in view of \eqref{eq:6.4} is equivalent to
\begin{equation}\label{eq:6.6}
    1-2\mbox{sn}^2Kr=\alpha,\quad r\in Q
\end{equation}
with $r=\frac{m-n}{n}$. Since $w_0=-|w_0|$ one obtains
from~\eqref{eq:6.5} and the reduction formula (see~\cite[Table
XII]{Akh})
\begin{equation}\label{eq:6.7}
    1-2\mbox{sn}^2\frac{K'\ln w_0}{\pi}=1-\frac{2}{k^2\mbox{sn}^2(\frac{K'\ln |w_0|}{\pi}+iK')}
    =1-\frac{2}{k^2\mbox{sn}^2(\frac{K'\ln |w_0|}{\pi})}.
\end{equation}
Hence one obtains that $\omega(0)$ is a rational number $\frac{m}{n}$ if and only if
$\frac{K'}{\pi K}\ln |w_0|$ is a rational number, or, equivalently,
\begin{equation}\label{eq:6.8}
     1-\frac{2}{k^2\mbox{sn}^2(Kr')}=\frac{1}{\alpha}
\end{equation}
for $r'=\frac{m-n}{n}$.

These calculations lead to the following
\begin{prop}\label{prop:5.9}
Let a finite nonnegative measure  $\sigma$ on
$E=[-1,\alpha]\cup[\beta,1]$ be absolutely continuous
($d\sigma(t)=\rho(t)dt$) and satisfies the Szeg\"o
conditions~\eqref{Szego}.
Then:
\begin{enumerate}
\item[(i)] If $\alpha$ cannot be represented in the form~\eqref{eq:6.6}
for some $r\in Q$ then there is a subsequence of $[n/n]$ Pad\'e
approximants of ${\mathfrak F}$ which converges to ${\mathfrak F}$
locally uniformly on
$\dC\setminus([-1-\varepsilon,\alpha]\cup[\beta,1+\varepsilon])$
for some $\varepsilon>0$.

\item[(ii)] If $\alpha$ satisfies~\eqref{eq:6.6} for some
$r=\frac{m}{n}$ with $m,n\in\dN$ $(gcd(m,n)=1)$ and does not
satisfies~\eqref{eq:6.6} for any $r'=\frac{k}{n}$ with $k\in\dN$
$|k|\le n$, then the sequence $\{{\mathfrak
F}^{[n-1/n-1]}\}_{n\in\cN_F}$
of Pad\'e approximants of ${\mathfrak F}$
converges
to ${\mathfrak F}$ locally uniformly on
$\dC\setminus([-1-\varepsilon,\alpha]\cup[\beta,1+\varepsilon])$ for
some $\varepsilon>0$.
\end{enumerate}
\end{prop}


In this paper we considered the case of one "turning point". The
following example shows that in the case of 2 "turning points" the
behaviour of diagonal Pad\'e approximants seems to be more
complicated.
\begin{ex}[\cite{Stahl}]
Let $\theta_1$, $\theta_2$, $1$ ($0<\theta_1<\theta_2<1$) be
rationally independent real, and let
\[
F(\lambda)=\int_{-1}^{1}\frac{1}
{t-\lambda}\frac{(t-\cos\pi\theta_1)(t-\cos\pi\theta_2)}
{\sqrt{1-t^2}}dt.
\]
Then all the diagonal Pad\'e approximants $F^{[k/k]}$ exist, but do
not converge locally uniformly on $\dC\setminus\dR$ since
\[
\bigcap\limits_{n=1}^{\infty}\overline{ \bigcup\limits_{k\ge
n}\cP(F^{[k/k]})}=\dC.
\]
\end{ex}

\end{document}